\documentclass[a4paper, final, 12pt]{article}
\usepackage[utf8]{inputenc}
\usepackage{amstext,amssymb,amsmath,amsbsy}
\usepackage[pagewise]{lineno}
\usepackage{hyperref}
\usepackage{amscd}
\usepackage{amsfonts}
\usepackage{indentfirst}
\usepackage{verbatim}
\usepackage{amsmath}
\usepackage{amsthm}
\usepackage{enumerate}
\usepackage{graphicx}
\usepackage{color}
\usepackage[OT1]{fontenc}
\usepackage{subfig}
\usepackage{algorithm}
\usepackage{algpseudocode}

\title{Travel time tomography with formally determined incomplete data in 3D}
\author{ Michael V. Klibanov\thanks{%
Department of Mathematics and Statistics, University of North Carolina
Charlotte, Charlotte, NC, 28223, mklibanv@uncc.edu.} }
\date{}

\begin{document}

\maketitle

\begin{abstract}
For the first time, a globally convergent numerical method is developed and
Lipschitz stability estimate is obtained for the challenging problem of
travel time tomography in 3D for formally determined incomplete data. The
semidiscrete case is considered meaning that finite differences are involved
with respect to two out of three variables. First, Lipschitz stability
estimate is derived, which implies uniqueness. Next, a weighted globally
strictly convex Tikhonov-like functional is constructed using a
Carleman-like weight function for a Volterra integral operator. The gradient
projection method is constructed to minimize this functional. It is proven
that this method converges globally to the exact solution if the noise in
the data tends to zero.
\end{abstract}

%\articletype{ARTICLE TEMPLATE}% Specify the article type or omit as appropriate

%\begin{keywords}
\noindent \textit{Key words.} inverse kinematic problem; Carleman-like
estimate for the Volterra operator; Lipschitz stability; convexification;
globally convergent numerical method %\end{keywords}

%\begin{amscode}
\noindent \textit{AMS Classification} 35R30 %\end{amscode}

\section{Introduction}

\label{sec:1}

We construct a globally convergent numerical method for the challenging 3D
travel time tomography problem (TTTP) with formally determined incomplete
data. The TTTP was first considered by Herglotz \cite{Herg} in 1905 and then
by Wiechert and Zoeppritz \cite{W} in 1907 in the 1D case due to a
geophysical application, also, see a detailed description of the 1D case in 
\cite{R}. However, globally convergent numerical methods for the 3D TTTP
with formally determined data were not developed so far. By the definition
of, e.g. \cite{BakKlib,BK}, a numerical method for a nonlinear problem is
called \emph{globally convergent} if there exists a theorem claiming that
this method delivers at least one point in a sufficiently small neighborhood
of the correct solution without relying on any advanced knowledge of this
neighborhood, i.e. a good first guess for the solution is not required.\emph{%
\ }

In this paragraph, we indicate those ideas for the TTTP, which are presented
here \emph{for the first time}. More details about the latter statement,
including some references, are given in this section below. Since we develop
a numerical method, we are allowed to work here with an \emph{approximate
mathematical model. }We study the case when the data for the TTTP are both
formally determined and incomplete. The TTTP is considered in the
semidiscrete form, i.e. \ we consider the practically important case of
finite differences with respect to two out of three variables, also see
Remark 6.1 in section 6. The Lipschitz stability estimate is obtained, which
implies uniqueness. Next, the exact solution is constructed via introducing
a sequence, which converges to that solution starting from an arbitrary
point of a certain bounded set, provided that the noise in the data tends to
zero. Since no restrictions are imposed on the size $R$ of that set, then
this is a globally convergent numerical method. To construct that method, we
introduce a weighted Tikhonov-like functional for the TTTP, which is
strictly convex on that bounded set, i.e. we use the so-called
\textquotedblleft convexification" procedure. However, while in all previous
versions of the convexification the so-called Carleman Weight Function (CWF)
was applied to differential operators, in the current paper the CWF is
applied to a Volterra-like nonlinear integral operator. Our method does not
use sophisticated geometrical properties to construct the above sequence.
Rather, we straightforwardly minimize the above mentioned globally strictly
convex Tikhonov-like functional.

The TTTP is the problem of the recovery of the spatially distributed speed
of acoustic waves from first times of arrival of those waves. Another well
known term for the TTTP is the \textquotedblleft inverse kinematic problem
of seismic" \cite{R}. Waves are originated by some point sources located on
the boundary of the domain of interest. First times of arrival are recorded
on a number of detectors located on that boundary. It is well known that the
TTTP is actually a nonlinear problem of the integral geometry, see, e.g. 
\cite{R}\emph{. }The\emph{\ }TTTP has important applications in geophysics 
\cite{Herg,R,Vol,W}. In addition, it was established in \cite{KR} that the
TTTP arises in the phaseless inverse problem of scattering of
electromagnetic waves at high frequencies. The specific TTTP considered here
has potential applications in geophysics, checking the bulky baggage in
airports, search for possible defects inside the walls, etc..

The \textquotedblleft formally determined data" means that the number $m$ of
free variables in the data equals to the number $n$ of free variables in the
unknown coefficient, $m=n$. If, however, $m>n$, then the data are
overdetermined. In our case $n=m=3$. In the 2D case all previously known
results for the TTTP work with formally determined data with $m=n=2$. Only
complete data for the TTTP were considered in the past. In the 3D case those
data were overdetermined with $m=4>n=3$. Complete data for the TTTP are
generated by the point source running along the entire surface of the
boundary of the domain of interest. Unlike these, our point source runs
along an interval of a straight line located outside of the domain of
interest. This means that the data are incomplete. The sole purpose of
Figure 1 is to illustrate this for a simple case when geodesics are straight
lines.

\begin{figure}[tbp]
\begin{center}
\subfloat[]{\includegraphics[width=0.3\textwidth]{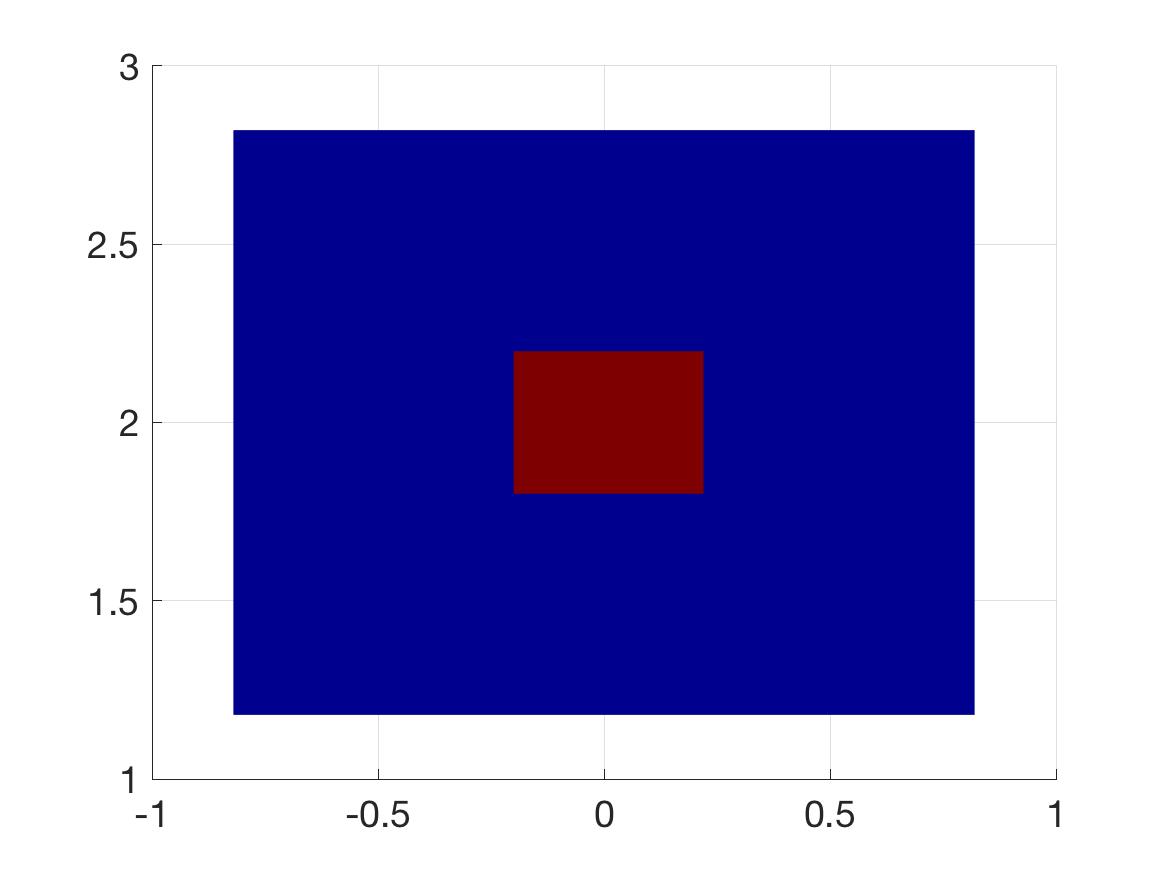}} \hfill %
\subfloat[]{\includegraphics[width=0.3\textwidth]{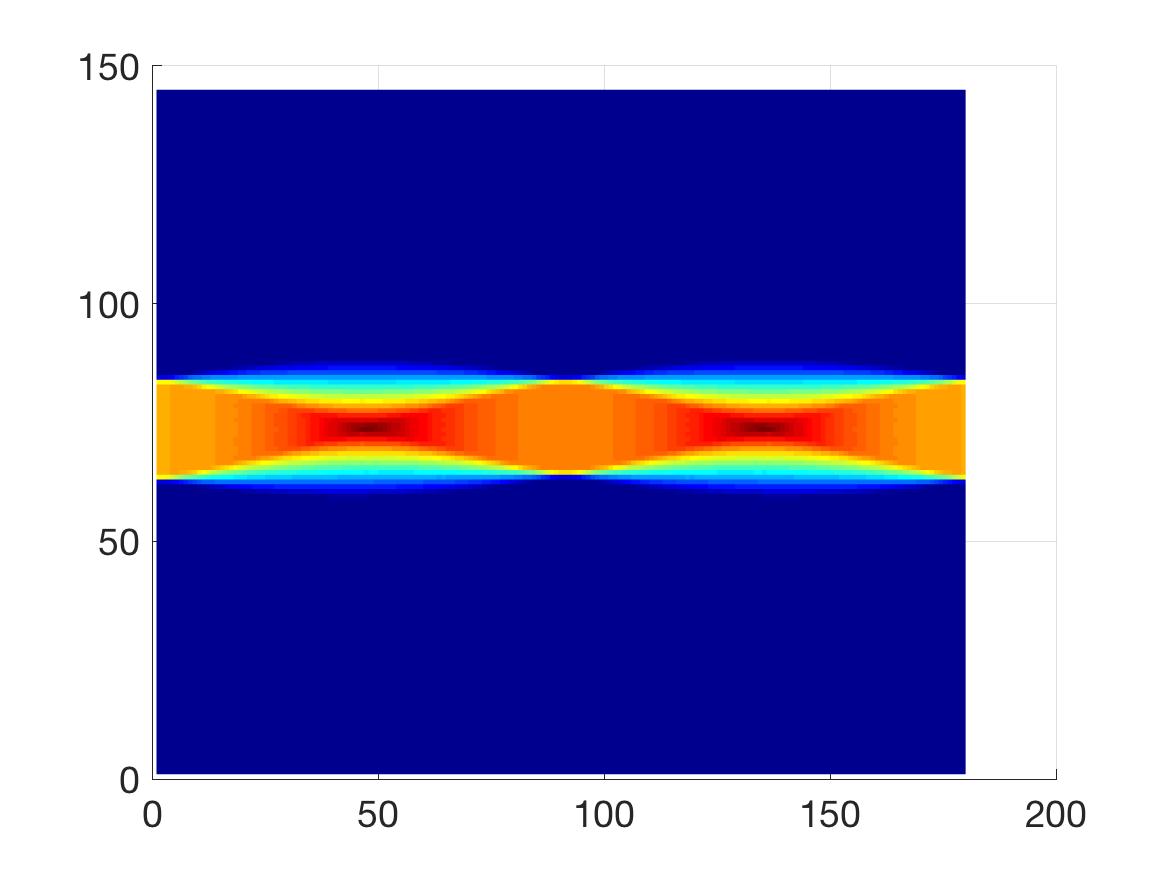}}
\hfill \subfloat[]{\includegraphics[width=0.3%
\textwidth]{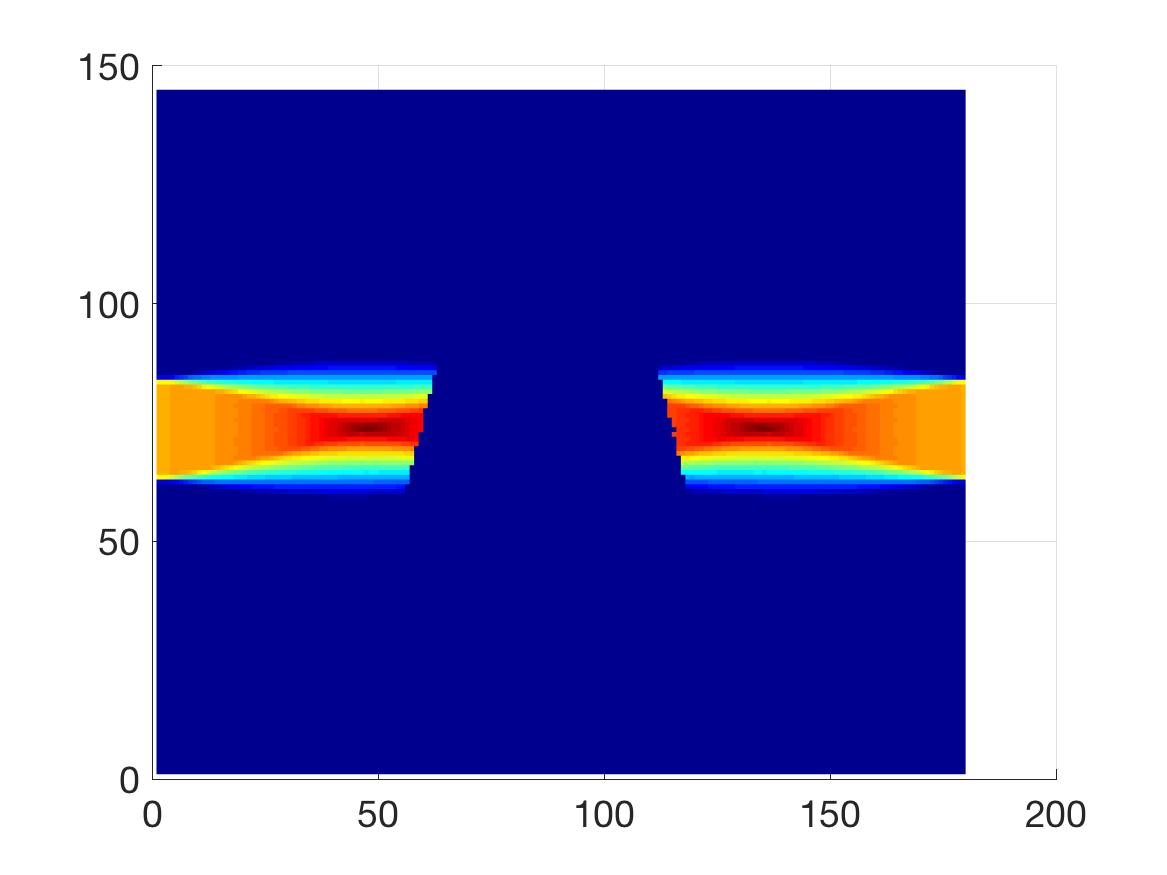}}
\end{center}
\caption{\emph{An illustration for complete and incomplete data in the 2D
case, see details in  \cite{R}. To simplify, we assume in this figure
that the geodesics are straight lines. Thus, we deal in this figure with the
data of Radon transform, generated by the function \textquotedblleft radon"
of MATLAB. a) The true function }$m\left( \mathbf{x}\right) $ \emph{to be
imaged. b) The complete data of the Radon transform of the function of a).
c) The incomplete data of the Radon transform of a) in the case when the
source runs along an interval of a straight line, as in this paper below.}}
\end{figure}

Our \emph{approximate mathematical model} consists of two items, see Remarks
6.2, 6.3 in the end of section 6 for more details. First, we consider a
semidiscrete model. This means that partial derivatives with respect to two
out of three variables are written in finite differences. The author
believes that this might be valuable for a possible future numerical
implementation of the method of this paper. Second, we assume that a certain
function generated by the solution of the eikonal equation can be
represented as a truncated Fourier series with respect to a special
orthonormal basis in the $L_{2}\left( 0,1\right) $ space. That basis was
first constructed in \cite{KJIIP}. Functions of that basis depend only on
the position of the source. The number $N\geq 1$ of terms of this series is
not allowed to tend to the infinity. We note that such assumptions quite
often take place in the theory of ill-posed problems, and corresponding
computational results are not affected by these assumptions, see, e.g. \cite%
{KEIT,convnew,KN}, where the same basis was used to obtain good quality
numerical results.

The TTTP is about the reconstruction of the right hand side of the nonlinear
eikonal Partial Differential Equation (PDE) from boundary measurements of
the solution of this equation. It is well known that conventional least
squares cost functionals for nonlinear inverse problems are non convex. The
non convexity, in turn causes the well known phenomenon of multiple local
minima and ravines of that functional, see, e.g. \cite{Scales} for a quite
convincing numerical example. Therefore, convergence of any minimization
algorithm to the correct solution is in question, unless its starting point
is chosen to be in a sufficiently small neighborhood of that solution, which
is the so-called small perturbation approach. However, it is unclear in many
applications how to actually obtain that sufficiently small neighborhood.

To avoid the phenomenon of local minima, we apply the so-called
convexification method. More precisely, we construct a weighted globally
strictly convex Tikhonov-like functional for our TTTP. The key element of
this functional is the presence in it of the weight function which looks
similar with the CWF.

The convexification was first proposed in \cite{Klib95,Klib97}.\ However,
those initial works lacked some important analytical results, which would
allow numerical studies. Such results were recently obtained in \cite%
{BakKlib}. Thus, recent publications \cite{KKN,KEIT,convnew} contain both
the theory and numerical results for the convexification method for some
coefficient inverse problems for PDEs, although not for the TTTP. We prove
the global convergence to the correct solution of the gradient projection
method of the minimization of our functional.

As to the numerical methods for the TTTP in the $n-$D case, $n=2,3$, such a
method for the 2D TTTP was published in \cite{Sch}. Another numerical
approach in 3D was published in \cite{Z}. Both publications \cite{Sch,Z}
use, at a certain stage, the minimization of a least squares cost
functional. Since the convexity of those functionals of \cite{Sch,Z} is not
proven, then the problem of local minima is not addressed there. In both
publications \cite{Sch,Z} complete data are used, and these data are
overdetermined ones in the 3D case of \cite{Z}.

The first global Lipschitz stability and uniqueness result for the TTTP was
obtained by Mukhometov in the 2D case \cite{Mukh}. Next, this result was
extended in \cite{BG,MR,R} to the $n-$D case, $n\geq 3$. We also refer to
the related work \cite{PU} for the 2D case. In all these references, the
data are complete and the assumption of the regularity of geodesic lines is
used. In addition, more recently the question of uniqueness in the 3D case
when geodesic lines are not necessarily regular ones was considered in \cite%
{Stef}. In the 2D case of \cite{Mukh,PU,Sch} the data are formally
determined. However, they are overdetermined in the $n-$D case with $n\geq 3$
\cite{BG,MR,R,Stef,Z}.

As to the Carleman estimates, for the first time they were introduced in the
field of Inverse Problems in \cite{BukhKlib} with the goal of proofs of
global uniqueness and stability results for coefficient inverse problems.
This idea became quite popular since then with many publications of many
authors. To shorten the citation list, we refer here only to, e.g. books 
\cite{BK,BY,KT}, the survey \cite{Ksurvey} and references cited therein. In
addition to uniqueness and stability, various modifications of the idea of 
\cite{BukhKlib} are currently used for the convexification, see
corresponding comments and references above.

All functions considered below are real valued ones. In section 2 we state
the TTTP. In section 3 we introduce a special orthonormal basis. In section
4 we estimate from the below a derivative of the solution of the eikonal
equation. In section 5 we derive a boundary value problem for a system of
coupled nonlinear integral differential equations. In section 6 we rewrite
that system in a semidiscrete form. In section 7, \ we establish Lipschitz
stability and uniqueness. In section 8 we construct the above mentioned
globally strictly convex functional and formulate corresponding theorems.
These theorems are proven in section 9.

\section{Statement of the problem}

\label{sec:2}

Below $\mathbf{x}=\left( x,y,z\right) \in \mathbb{R}^{3}.$ Let numbers $%
A,\sigma =const.>0$. Define the rectangular prism $\Omega \subset \mathbb{R}%
^{3}$ as%
\begin{equation}
\Omega =\left\{ \mathbf{x}=\left( x,y,z\right) :x,y\in \left( 0,1\right)
,z\in \left( A,A+\sigma \right) \right\} .  \label{1.0}
\end{equation}%
Denote parts of the boundary $\partial \Omega $ as%
\begin{equation}
B_{A}=\left\{ \mathbf{x}=\left( x,y,z\right) :x,y\in \left( 0,1\right)
,z=A\right\} ,  \label{1.1}
\end{equation}%
\begin{equation}
B_{A+\sigma }=\left\{ \mathbf{x}=\left( x,y,z\right) :x,y\in \left(
0,1\right) ,z=A+\sigma \right\} ,  \label{1.2}
\end{equation}%
\begin{equation}
\Gamma =\partial \Omega \diagdown \left( B_{A}\cup B_{A+\sigma }\right) .
\label{1.3}
\end{equation}

Let $n\left( \mathbf{x}\right) $ be the refractive index of the medium at
the point $\mathbf{x.}$ Hence, $c\left( \mathbf{x}\right) =1/n\left( \mathbf{%
x}\right) $ is the sound speed. Denote $m\left( \mathbf{x}\right)
=n^{2}\left( \mathbf{x}\right) .$ Let the number $m_{0}>0$ be given. We
impose the following assumptions on the function $m\left( \mathbf{x}\right)
: $%
\begin{equation}
m\left( \mathbf{x}\right) \geq m_{0},\quad \mathbf{x}\in \mathbb{R}^{3},
\label{2.2}
\end{equation}%
\begin{equation}
m\left( \mathbf{x}\right) =1,\quad \mathbf{x}\in \left\{ z<A\right\} ,
\label{2.3}
\end{equation}%
\begin{equation}
m\in C^{2}\left( \mathbb{R}^{3}\right) ,  \label{2.4}
\end{equation}%
\begin{equation}
m_{z}\left( \mathbf{x}\right) \geq 0\text{, }\mathbf{x}\in \overline{\Omega }%
.  \label{2.5}
\end{equation}

\textbf{Remark 2.1}. \emph{We note that the monotonicity condition (\ref{2.5}%
) is not an overly restrictive one.\ Indeed, it can also be found in the end
of section 2 of chapter 3 of the book \cite{R}: see formulas (3.24) and (3.24%
}$^{\prime }$\emph{) there. Also, an analogous monotonicity condition was
actually imposed in the 1D case in the originating classical works of
Herglotz and Wiechert and Zoeppritz \cite{Herg,W}, see section 3 of chapter
3 of \cite{R} for a description of their method. We also refer to Figures 5
and 10 in \cite{Vol} for some geophysical information. }

The function $m\left( \mathbf{x}\right) $ generates the Riemannian metric 
\begin{equation}
d\tau =\sqrt{m(\mathbf{x})}\left\vert d\mathbf{x}\right\vert ,|d\mathbf{x}|=%
\sqrt{(dx)^{2}+(dy)^{2}+(dz)^{2}}.  \label{2.6}
\end{equation}%
The travel time from the point $\mathbf{x}_{0}\in \mathbb{R}^{3}$ (source)
to the point $\mathbf{x}\in \mathbb{R}^{3}$ (receiver) is \cite{R}%
\begin{equation}
\tau \left( \mathbf{x},\mathbf{x}_{0}\right) =\int\displaylimits_{\Gamma \left( 
\mathbf{x},\mathbf{x}_{0}\right) }\sqrt{m\left( \mathbf{x}\left( s\right)
\right) }{ds},  \label{2.7}
\end{equation}%
where $\Gamma \left( \mathbf{x},\mathbf{x_{0}}\right) $ is the geodesic line
connecting points $\mathbf{x}$ and $\mathbf{x}_{0}$ and $ds$ is the
euclidean arc length. We assume that the source $\mathbf{x}_{0}$ runs along
an interval $L$ of a straight line located in the plane $\left\{ z=0\right\}
,$ 
\begin{equation}
L=\left\{ \mathbf{x=}\left( x,y,z\right) :x=\alpha \in \left( 0,1\right)
,y=1/2,z=0\right\} .  \label{2.8}
\end{equation}%
Hence, $\mathbf{x}_{0}=\mathbf{x}_{\alpha }=\left( \alpha ,1/2,0\right)
,\alpha \in \left( 0,1\right) .$ Let $\tau \left( \mathbf{x,}\alpha \right) $
be the travel time between points $\mathbf{x}$ and $\mathbf{x}_{\alpha
}=\left( \alpha ,1/2,0\right) .$ Thus, we denote $\tau \left( \mathbf{x,}%
\alpha \right) =\tau \left( \mathbf{x},\mathbf{x}_{\alpha }\right) .$ We
denote $\Gamma \left( \mathbf{x},\alpha \right) $ the geodesic line
connecting points $\mathbf{x}$ and $\mathbf{x}_{\alpha }.$ It is well known 
\cite{R} that the function $\tau \left( \mathbf{x},\alpha \right) $
satisfies eikonal equation as the function of $\mathbf{x}$, 
\begin{equation}
\tau _{x}^{2}+\tau _{y}^{2}+\tau _{z}^{2}=m\left( \mathbf{x}\right) ,\text{ }%
\mathbf{x\in }\mathbb{R}^{3},  \label{2.9}
\end{equation}%
\begin{equation*}
\tau \left( \mathbf{x},\alpha \right) =O\left( \left\vert \mathbf{x}-\mathbf{%
x}_{\alpha }\right\vert \right) ,\text{ }{\text{a}}\text{{s }}\left\vert 
\mathbf{x}-\mathbf{x}_{\alpha }\right\vert \rightarrow 0.
\end{equation*}%
Everywhere below we assume without further mentioning that the following
property holds:

\textbf{Regularity of Geodesic Lines}. \emph{The function }$\tau \left( 
\mathbf{x},\alpha \right) \in C^{2}\left( \mathbb{R}^{3}\times \left[ 0,1%
\right] \right) .$ \emph{For any pair of points }$\left( \mathbf{x},\mathbf{x%
}_{\alpha }\right) $\emph{\ }$\in \overline{\Omega }\times L$ \emph{there
exists unique geodesic line }$\Gamma \left( \mathbf{x},\alpha \right) $\emph{%
\ connecting these two points and }$\Gamma \left( \mathbf{x},\alpha \right)
\cap B_{A}\neq \varnothing $\emph{.} \emph{In addition, if any geodesic
line, which starts at a point }$\mathbf{x}_{\alpha }\in L,$ \emph{intersects 
}$B_{A},$ \emph{then it intersects it at a single point. Also, it does not
go \textquotedblleft downwards" in the }$z-$\emph{direction, but instead
intersects }$\partial \Omega \diagdown B_{A}$\emph{\ at another single
point, see (\ref{1.0})-(\ref{1.3}). In addition, after intersecting }$%
\partial \Omega \diagdown B_{A},$\emph{\ this line does not
\textquotedblleft come back" in the domain }$\Omega $\emph{\ but rather goes
away from this domain}$.$\emph{\ In other words, this line is not reflected
back from any point of its intersection with }$\partial \Omega $.

The following sufficient condition of the regularity of geodesic lines was
derived in \cite{R1}:%
\begin{equation*}
\sum_{i,j=1}^{3}\frac{\partial ^{2}\ln n(\mathbf{x})}{\partial x_{i}\partial
x_{j}}\xi _{i}\xi _{j}\geq 0,\>\forall \xi \in \mathbb{R}^{3},\forall 
\mathbf{x=}\left( x_{1},x_{2},x_{3}\right) \in \overline{\Omega }.
\end{equation*}

\textbf{Travel Time Tomography Problem (TTTP)}. \emph{Suppose that the
function }$m\left( \mathbf{x}\right) $ \emph{satisfies conditions (\ref{2.2}%
)-(\ref{2.5}). Assume that} \emph{the following function }$f\left( \mathbf{x}%
,a\right) $\emph{\ is given:}%
\begin{equation}
\tau \left( \mathbf{x,}\alpha \right) =f\left( \mathbf{x},\alpha \right)
,\forall \mathbf{x}\in \Gamma \cup B_{A+\sigma },\forall \alpha \in \left(
0,1\right) .  \label{2.10}
\end{equation}%
\emph{Determine the function} $m\left( \mathbf{x}\right) $ \emph{for} $%
\mathbf{x}\in \Omega .$

In other words, by (\ref{1.0})-(\ref{1.3}) and (\ref{2.10}) the data for the
travel time are given for all sources running along the line interval $L$
defined in (\ref{2.8}) and for the part $\Gamma \cup B_{A+\sigma }$ of the
boundary $\partial \Omega .$ Hence, the data (\ref{2.10}) are both formally
determined and incomplete.

\section{A Special Orthonormal Basis}

\label{sec:3}

We now reproduce a special orthonormal basis $\left\{ \Psi _{n}\left( \alpha
\right) \right\} _{n=0}^{\infty }$ in $L_{2}\left( 0,1\right) ,$ which was
constructed in \cite{KJIIP}. This basis has the following two properties:

\begin{enumerate}
\item $\Psi _{n}\in C^{1}\left[ 0,1\right] ,$ $\forall n=1,2,\dots $

\item Let $\left\{ ,\right\} $ be the scalar product in $L_{2}\left(
0,1\right) $. Denote $a_{mn}=\left\{ \Psi _{n}^{\prime },\Psi _{m}\right\} .$
Then the matrix $M_{N}=\left( a_{mn}\right) _{m,n=0}^{N}$ should be
invertible for any $N=1,2,\dots $
\end{enumerate}

Note that if one would use either the basis of standard orthonormal
polynomials or the trigonometric basis, then the first derivative of its
first element would be identically zero. Hence, the matrix $M_{N}$ would not
be invertible in this case.

We now describe the basis of \cite{KJIIP}. Consider the set of functions $%
\left\{ \xi _{n}\left( \alpha \right) \right\} _{n=0}^{\infty }=\{\left(
\alpha +a\right) ^{n}e^{\alpha }\}_{n=0}^{\infty }$, where $a=const.>0.$
This is a set of linearly independent functions. Besides, this set is
complete in the space $L_{2}(0,1)$. After applying the Gram-Schmidt
orthonormalization procedure to this set, we obtain the orthonormal basis $%
\{\Psi _{n}\left( \alpha \right) \}_{n=1}^{\infty }$ in $L_{2}(0,1)$. In
fact, the function $\Psi _{n}(\alpha )$ has the form $\Psi _{n}(\alpha
)=P_{n}(\alpha +a)e^{\alpha }$, $\forall n\geq 0,$ where $P_{n}$ is a
polynomial of the degree $n$. Thus, functions $\Psi _{n}(\alpha )$ are
polynomials orthonormal in $L_{2}(0,1)$ with the weight $e^{2\alpha }.$ The
matrix $M_{N}$ is invertible since its elements $a_{mn}=\left( \Psi
_{n}^{\prime },\Psi _{m}\right) $ are such that $a_{mn}=1$ if $m=n$ and $%
a_{mn}=0$ if $n<m$. \newline

Consider the function $q\left( \alpha \right) $ in the following form:%
\begin{equation}
q\left( \alpha \right) =\sum\displaylimits_{n=1}^{N}q_{n}\Psi _{n}(\alpha ),\text{ 
}q_{n}=\int\displaylimits_{0}^{1}q\left( \alpha \right) \Psi _{n}(\alpha )d\alpha .
\label{1}
\end{equation}%
Below we need to impose such a sufficient condition on the vector of
coefficients

$q^{N}=\left( q_{0},...,q_{N-1}\right) ^{T}$ in the Fourier expansion (\ref%
{1}) which would guarantee that the function $q\left( \alpha \right) $ is
positive for all $\alpha \in \left[ 0,1\right] .$ Consider vector functions $%
\xi ^{N}\left( \alpha \right) =\left( \xi _{1},...,\xi _{N}\right)
^{T}\left( \alpha \right) $ and

$\Psi ^{N}(\alpha )=\left( \Psi _{1},...,\Psi _{N-1}\right) ^{T}\left(
\alpha \right) .$ The desired condition is given in Lemma 3.1.

\textbf{Lemma 3.1.} \emph{Let the }$N\times N$\emph{\ matrix }$X_{N}$\emph{\
transforms the vector function }$\xi ^{N}\left( \alpha \right) $\emph{\ in
the vector function }$\Psi ^{N}\left( \alpha \right) $\emph{\ via the
Gram-Schmidt orthonormalization procedure, i.e. }$X_{N}\cdot \xi ^{N}\left(
\alpha \right) =\Psi ^{N}\left( \alpha \right) .$\emph{\ Let the matrix }$%
X_{N}^{T}$\emph{\ be the transpose of }$X_{N}$\emph{. Let }$X_{N}^{T}q^{N}=%
\widetilde{q}^{N}=\left( \widetilde{q}_{0},...,\widetilde{q}_{N-1}\right)
^{T}.$\emph{\ Suppose that all numbers }$\widetilde{q}_{0},...,\widetilde{q}%
_{N-1}>0.$\emph{\ Then in (\ref{1}) the function }$q\left( \alpha \right) >0$%
\emph{\ for all }$\alpha \in \left[ 0,1\right] .$

\textbf{Proof.} It follows from the Gram-Schmidt procedure that elements of
the matrix $X_{N}$ are independent on $\alpha .$ Let the raw number $n$ of
the matrix $X_{N}$ be $\left( x_{n1},x_{n2},...,x_{nN}\right) ,n=1,...,N.$
Then%
\begin{equation*}
\Psi _{n}\left( \alpha \right) =\sum\displaylimits_{j=1}^{N}x_{nj}\xi _{j}\left(
\alpha \right) .
\end{equation*}%
Hence, by (\ref{1}) 
\begin{equation}
q\left( \alpha \right)
=\sum\displaylimits_{n=1}^{N}q_{n}\sum\displaylimits_{j=1}^{N}x_{nj}\xi _{j}\left(
\alpha \right) =\sum\displaylimits_{j=1}^{N}\left(
\sum\displaylimits_{n=1}^{N}x_{nj}q_{n}\right) \xi _{j}\left( \alpha \right)
=\sum\displaylimits_{j=1}^{N}\widetilde{q}_{j}\xi _{j}\left( \alpha \right) .
\label{2}
\end{equation}%
Since $\xi _{j}\left( \alpha \right) =\left( \alpha +a\right) ^{j}e^{\alpha
}>0,$ then (\ref{2}) implies that $q\left( \alpha \right) >0$ for $\alpha
\in \left[ 0,1\right] .$ $\square $

\section{Estimate of $ \tau _{z}^{2}\left( \mathbf{x},\mathbf{ \alpha }\right) $ From the Below}

\label{sec:4}

\textbf{Lemma 4.1}. \emph{Assume that conditions (\ref{2.2})-(\ref{2.5})
hold.} \emph{Then} 
\begin{equation}
\tau _{z}^{2}\left( \mathbf{x},\alpha \right) \geq \frac{A^{2}}{A^{2}+2},%
\text{ }\forall \mathbf{x}\in \overline{\Omega },\forall \alpha \in \left[
0,1\right] .  \label{3.01}
\end{equation}%
\emph{Also,} 
\begin{equation}
\tau _{z}\left( \mathbf{x},\alpha \right) >0,\text{ }\forall \mathbf{x}\in 
\overline{\Omega }\cup \left\{ z\in \left( 0,A\right] \right\} ,\forall
\alpha \in \left[ 0,1\right] .  \label{300}
\end{equation}%
\emph{Thus,}%
\begin{equation}
\tau _{z}\left( \mathbf{x},\alpha \right) \geq \frac{A}{\sqrt{A^{2}+2}}%
,\forall \mathbf{x}\in \overline{\Omega },\forall \alpha \in \left[ 0,1%
\right] .  \label{301}
\end{equation}

\textbf{Proof}. Note that having proven (\ref{3.01}) is not enough for our
technique: we need to know the sign of the function $\tau _{z}\left( \mathbf{%
x},\alpha \right) ,$ i.e. (\ref{300}), in section 5 (more precisely, in (\ref%
{4.01})) where we consider the square root of\emph{\ }$\tau _{z}^{2}\left( 
\mathbf{x},\alpha \right) .$ Denote%
\begin{equation}
p=\tau _{x}(x,y,z,\alpha ),q=\tau _{x}(x,y,z,\alpha ),r=\tau
_{z}(x,y,z,\alpha ).  \label{3.1b}
\end{equation}%
The following equations for geodesic lines can be found on page 66 of \cite%
{R}:%
\begin{equation}
\frac{dx}{ds}=\frac{p}{m},\quad \frac{dy}{ds}=\frac{q}{m},\quad \frac{dz}{ds}%
=\frac{r}{m},  \label{3.2a}
\end{equation}%
\begin{equation}
\frac{dp}{ds}=\frac{m_{x}}{2m},\quad \frac{dq}{ds}=\frac{m_{y}}{2m},\quad 
\frac{dr}{ds}=\frac{m_{z}}{2m},  \label{3.2b}
\end{equation}%
where $s$ is a parameter. Using (\ref{3.1b}), we obtain\emph{\ }for $\tau
=\tau \left( x\left( s\right) ,y\left( s\right) ,z\left( s\right) ,\alpha
\right) $ along a geodesic line 
\begin{equation}
\frac{d\tau }{ds}=\frac{\partial \tau }{\partial x}\frac{dx}{ds}+\frac{%
\partial \tau }{\partial y}\frac{dy}{ds}+\frac{\partial \tau }{\partial z}%
\frac{dz}{ds}=p\frac{p}{m}+q\frac{q}{m}+r\frac{r}{m}=1.  \label{3.2c}
\end{equation}%
Set 
\begin{equation}
\tau \left( x\left( 0\right) ,y\left( 0\right) ,z\left( 0\right) ,\alpha
\right) =0\text{ for }s=0.  \label{3.2d}
\end{equation}
Then (\ref{3.2c}) implies:

$\tau \left( x\left( s\right) ,y\left( s\right) ,z\left( s\right) ,\alpha
\right) =s$. Hence, the parameter $s$ coincides with the travel time. In
particular, we now rewrite equations (\ref{3.2a}), (\ref{3.2b}) in a
different form: to have derivatives with respect to $z$ rather than with
respect to $s$. Hence, we obtain from (\ref{3.2a}) and (\ref{3.2b}): 
\begin{equation}
\frac{dx}{dz}=\frac{p}{r},\quad \frac{dy}{dz}=\frac{q}{r},\quad \frac{dp}{dz}%
=\frac{m_{x}}{2r},\quad \frac{dq}{dz}=\frac{m_{y}}{2r},\quad \frac{dr^{2}}{dz%
}=m_{z},\quad \frac{d\tau }{dz}=\frac{m}{r}.  \label{3.2}
\end{equation}%
\begin{equation}
x|_{z=0}=\alpha ,\quad y|_{z=0}=1/2,\quad p|_{z=0}=p_{0},\quad
q|_{z=0}=q_{0},\quad r|_{z=0}=r_{0},\text{ \ }\tau \mid _{z=0}=0.
\label{3.16}
\end{equation}%
where $\alpha \in \left( 0,1\right) $ and $p_{0}$, $q_{0}$ are some given
numbers such that $p_{0}^{2}+q_{0}^{2}\leq 1.$ The latter inequality follows
from (\ref{2.9}) and the fact that by (\ref{2.3}) $m\left( x,y,0\right) =1$.
Also, by (\ref{2.9}) $r_{0}=\pm \sqrt{1-p_{0}^{2}-q_{0}^{2}}.$ To prove that
we should take \textquotedblleft +" sign in the latter formula, we note that 
\begin{equation}
\tau \left( x,y,z,\alpha \right) =\sqrt{\left( x-\alpha \right) ^{2}+\left(
y-1/2\right) ^{2}+z^{2}}\text{ for }\left( x,y,z\right) \in \left(
0,1\right) \times \left( 0,1\right) \times \left( 0,A\right) .  \label{3.161}
\end{equation}%
Hence, $\tau _{z}=r=z/\tau >0$ for $z\in \left( 0,A\right) .$ Hence, 
\begin{equation}
r_{0}=\sqrt{1-p_{0}^{2}-q_{0}^{2}}\geq 0,  \label{3.160}
\end{equation}%
\begin{equation}
\tau _{z}^{2}(x,y,A,\alpha )\geq \frac{A^{2}}{A^{2}+2}\text{ for }\left(
x,y\right) \in \left( 0,1\right) \times \left( 0,1\right) ,\alpha \in \left(
0,1\right) .  \label{3.17}
\end{equation}

Suppose that the geodesic line defined by (\ref{3.2}) and (\ref{3.16})
intersects the part $B_{A}$ of the boundary $\partial \Omega .$ Then the
condition of the regularity of geodesic lines implies that there exists a
single number $z_{0}=z_{0}\left( p_{0},q_{0},\alpha \right) \in \left(
A,A+\sigma \right] $ such that the point $\left( x\left( z_{0}\right)
,y\left( z_{0}\right) ,z_{0}\right) \in \partial \Omega \diagdown B_{A}$ and
for all numbers $z\in \left( A,z_{0}\right) $ all points $\left( x\left(
z\right) ,y\left( z\right) ,z\right) $ of that geodesic line belong to $%
\Omega $. Since by (\ref{3.2}) $dr^{2}/dz=m_{z},$ then, using (\ref{3.2})
and (\ref{3.17}), we obtain 
\begin{equation*}
r^{2}\left( x\left( z\right) ,y\left( z\right) ,z,\alpha \right)
=\int\displaylimits_{A}^{z}m_{z}\left( x\left( t\right) ,y\left( t\right)
,t,\alpha \right) dt+r^{2}\left( x\left( A\right) ,y\left( A\right) ,A\right)
\end{equation*}%
\begin{equation*}
\geq r^{2}\left( x\left( A\right) ,y\left( A\right) ,A\right) \geq \frac{%
A^{2}}{A^{2}+2},z\in \left( A,z_{0}\right) ,\text{ }
\end{equation*}%
which establishes (\ref{3.01}). To prove (\ref{300}), we notice that it
follows from (\ref{2.2}), (\ref{2.5}), (\ref{3.1b}), the last equation (\ref%
{3.2b}) and (\ref{3.2d}) that 
\begin{equation}
\tau _{z}\left( \left( x\left( s\right) ,y\left( s\right) ,z\left( s\right)
,\alpha \right) \right) =\int\displaylimits_{0}^{s}\left( \frac{m_{z}}{2m}\right)
\left( x\left( t\right) ,y\left( t\right) ,z\left( t\right) ,\alpha \right)
dt\geq 0.  \label{302}
\end{equation}%
Estimates (\ref{300}) and (\ref{301}) follow immediately from (\ref{3.01}), (%
\ref{3.161}) and (\ref{302}). $\square $

\section{A Boundary Value Problem for a System of Nonlinear Coupled Integro
Differential Equations}

\label{sec:5}

\subsection{A nonlinear integro differential equation}

\label{sec:5.1}

Denote 
\begin{equation}
u\left( \mathbf{x},\alpha \right) =\tau _{z}^{2}\left( \mathbf{x},\alpha
\right) ,\mathbf{x}\in \Omega ,\alpha \in \left( 0,1\right) .  \label{4.1}
\end{equation}%
By (\ref{300}) 
\begin{equation}
\tau _{z}\left( \mathbf{x},\alpha \right) =\sqrt{u\left( \mathbf{x},\alpha
\right) }\text{, }\mathbf{x}\in \Omega .  \label{4.01}
\end{equation}%
Hence, (\ref{2.10}) and (\ref{4.01}) imply that for all $\alpha \in \left(
0,1\right) $ 
\begin{equation}
\tau \left( x,y,z,\alpha \right) =-\int\displaylimits_{z}^{A+\sigma }\sqrt{u\left(
x,y,t,\alpha \right) }dt+f\left( x,y,A+\sigma ,\alpha \right) ,\left(
x,y,z\right) \in \Omega ,  \label{4.2}
\end{equation}%
\begin{equation}
\tau _{x}\left( x,y,z,\alpha \right) =-\int\displaylimits_{z}^{A+\sigma }\left( 
\frac{u_{x}}{2\sqrt{u}}\right) \left( x,y,t,\alpha \right) dt+f_{x}\left(
x,y,A+\sigma ,\alpha \right) ,\left( x,y,z\right) \in \Omega ,  \label{4.3}
\end{equation}%
\begin{equation}
\tau _{y}\left( x,y,z,\alpha \right) =-\int\displaylimits_{z}^{A+\sigma }\left( 
\frac{u_{y}}{2\sqrt{u}}\right) \left( x,y,t,\alpha \right) dt+f_{y}\left(
x,y,A+\sigma ,\alpha \right) ,\left( x,y,z\right) \in \Omega .  \label{4.4}
\end{equation}%
Substituting (\ref{4.1})-(\ref{4.4}) in the eikonal equation (\ref{2.9}), we
obtain the following equation for $\left( x,y,z\right) \in \Omega ,\alpha
\in \left( 0,1\right) $: 
\begin{equation*}
u\left( x,y,z,\alpha \right) +\left[ -\int\displaylimits_{z}^{A+\sigma }\left( 
\frac{u_{x}}{2\sqrt{u}}\right) \left( x,y,t,\alpha \right) dt+f_{x}\left(
x,y,A+\sigma ,\alpha \right) \right] ^{2}
\end{equation*}%
\begin{equation}
+\left[ -\int\displaylimits_{z}^{A_{1}+\sigma }\left( \frac{u_{y}}{2\sqrt{u}}%
\right) \left( x,y,t,\alpha \right) dt+f_{y}\left( x,y,A+\sigma ,\alpha
\right) \right] ^{2}=m\left( x,y,z\right)  \label{4.5}
\end{equation}%
Differentiating (\ref{4.5}) with respect to $\alpha $ and using $\partial
_{\alpha }m\left( x,y,z\right) \equiv 0,$ we obtain for $\left( x,y,z\right)
\in \Omega ,\alpha \in \left( 0,1\right) $ 
\begin{equation*}
u_{\alpha }\left( x,y,z,\alpha \right) +\frac{\partial }{\partial \alpha }%
\left[ -\int\displaylimits_{z}^{A+\sigma }\left( \frac{u_{x}}{2\sqrt{u}}\right)
\left( x,y,t,\alpha \right) dt+f_{x}\left( x,y,A+\sigma ,\alpha \right) %
\right] ^{2}
\end{equation*}%
\begin{equation}
+\frac{\partial }{\partial \alpha }\left[ -\int\displaylimits_{z}^{A+\sigma
}\left( \frac{u_{y}}{2\sqrt{u}}\right) \left( x,y,t,\alpha \right)
dt+f_{x}\left( x,y,A+\sigma ,\alpha \right) \right] ^{2}.  \label{4.6}
\end{equation}

\subsection{Boundary value problem for a system of coupled integro
differential equations}

\label{sec:5.2}

Using (\ref{3.161}) and (\ref{4.1}), denote for $\left( x,y,\alpha \right)
\in \left[ 0,1\right] ^{3}$ 
\begin{equation}
u_{0}\left( \mathbf{x},\alpha \right) =u(x,y,A,\alpha )=\frac{A^{2}}{\left(
x-\alpha \right) ^{2}+\left( y-1/2\right) ^{2}+A^{2}}\geq \frac{A^{2}}{%
A^{2}+2}.  \label{4.8}
\end{equation}%
We seek the function $u\left( \mathbf{x},\alpha \right) $ in the form%
\begin{equation}
u\left( \mathbf{x},\alpha \right) =u_{0}\left( \mathbf{x},\alpha \right)
+v\left( \mathbf{x},\alpha \right) ,\mathbf{x}\in \overline{\Omega },\alpha
\in \left[ 0,1\right] ,  \label{4.9}
\end{equation}%
\begin{equation}
v\left( x,y,A,\alpha \right) =0,\left( x,y,\alpha \right) \in \left[ 0,1%
\right] ^{3},  \label{4.10}
\end{equation}%
where the function $v\left( \mathbf{x},\alpha \right) $ is unknown for $%
\mathbf{x}\in \Omega ,\alpha \in \left( 0,1\right) .$ Recall that the part $%
\Gamma $ of the boundary of the domain $\Omega $ is defined in (\ref{1.3}).
We need to obtain zero boundary condition at $\Gamma $ for a function
associated with the function $v$. To do this, we assume first that there
exists a function $g\left( \mathbf{x},\alpha \right) \in H^{1}\left( \Omega
\right) $ for every $\alpha \in \left[ 0,1\right] $ such that 
\begin{equation}
g\left( \mathbf{x},\alpha \right) =\left( f_{z}\right) ^{2}\left( \mathbf{x}%
,\alpha \right) -u_{0}\left( \mathbf{x},\alpha \right) ,\forall \mathbf{x\in 
}\Gamma \mathbf{,\forall }\alpha \in \left[ 0,1\right] ,  \label{4.11}
\end{equation}%
\begin{equation}
g\left( x,y,A,\alpha \right) =0.  \label{4.110}
\end{equation}%
Introduce the function $w\left( \mathbf{x},\alpha \right) ,$ 
\begin{equation}
w\left( \mathbf{x},\alpha \right) =v\left( \mathbf{x},\alpha \right)
-g\left( \mathbf{x},\alpha \right) ,\mathbf{x}\in \overline{\Omega },\alpha
\in \left[ 0,1\right] .  \label{4.12}
\end{equation}%
Then (\ref{4.9})-(\ref{4.12}) imply that 
\begin{equation}
u\left( \mathbf{x},\alpha \right) =u_{0}\left( \mathbf{x},\alpha \right)
+w\left( \mathbf{x},\alpha \right) +g\left( \mathbf{x},\alpha \right) ,%
\mathbf{x}\in \overline{\Omega },\alpha \in \left[ 0,1\right] ,  \label{4.13}
\end{equation}%
\begin{equation}
w\left( \mathbf{x},\alpha \right) =0,\forall \mathbf{x\in }\Gamma \mathbf{%
,\forall }\alpha \in \left[ 0,1\right] ,  \label{4.14}
\end{equation}%
\begin{equation}
w\left( x,y,A,\alpha \right) =0,\left( x,y,\alpha \right) \in \left[ 0,1%
\right] ^{3}.  \label{4.15}
\end{equation}

We assume that both functions $w\left( \mathbf{x},\alpha \right) $ and $%
g\left( \mathbf{x},\alpha \right) $ have the form of the truncated Fourier
series with respect to the orthonormal basis $\{\Psi _{n}\left( \alpha
\right) \},$%
\begin{equation}
w\left( \mathbf{x},\alpha \right) =\sum\displaylimits_{n=1}^{N}w_{n}\left( \mathbf{%
x}\right) \Psi _{n}\left( \alpha \right) ,\mathbf{x}\in \overline{\Omega }%
,\alpha \in \left[ 0,1\right] ,  \label{4.16}
\end{equation}%
\begin{equation}
g\left( \mathbf{x},\alpha \right) =\sum\displaylimits_{n=1}^{N}g_{n}\left( \mathbf{%
x}\right) \Psi _{n}\left( \alpha \right) ,\mathbf{x}\in \overline{\Omega }%
,\alpha \in \left[ 0,1\right] .  \label{4.17}
\end{equation}%
Here, coefficients $w_{n}\left( \mathbf{x}\right) $ are unknown and
coefficients $g_{n}\left( \mathbf{x}\right) $ are known. And similarly for $%
w_{x},w_{y},w_{x\alpha },w_{y\alpha }$ and the same derivatives of the
function $g$. Furthermore, we assume that these functions, being substituted
in equation (\ref{4.6}), give us zero in its right hand side for $\mathbf{x}%
\in \Omega ,\alpha \in \left( 0,1\right) .$\emph{\ }By (\ref{4.14})-(\ref%
{4.16}) 
\begin{equation}
w_{n}\left( \mathbf{x}\right) \mid _{\Gamma }=0,w_{n}\left( x,y,A\right) =0.
\label{4.18}
\end{equation}%
Denote

\begin{equation}
W\left( \mathbf{x}\right) =\left( w_{1},...,w_{N}\right) ^{T}\left( \mathbf{x%
}\right) \text{,}  \label{4.19}
\end{equation}%
\begin{equation}
G\left( \mathbf{x}\right) =\left( g_{1},...,g_{N}\right) ^{T}\left( \mathbf{x%
}\right) ,  \label{4.20}
\end{equation}

\textbf{Remark 5.1}. \emph{Note that we do not impose an analog of
conditions (\ref{4.16}), (\ref{4.17}) on the function} $u(x,y,A,\alpha ).$ 
\emph{The number }$N$\emph{\ of Fourier harmonics in (\ref{4.16}), (\ref%
{4.17}) should be chosen in numerical studies. For example, analogs of the
series (\ref{4.16}) were considered for four different inverse problems in 
\cite{KEIT,convnew,KN,SKN}. It was established numerically that for an
inverse problem of \cite{KEIT} the optimal }$N=8$\emph{, for the inverse
problem of \cite{convnew} the optimal }$N=3$\emph{, for the inverse problem
of \cite{KN} the optimal }$N=15$\emph{\ and the optimal }$N=12$\emph{\ for
the inverse problem of \cite{SKN}. }

Let 
\begin{equation}
f\left( x,y,A+\sigma ,\alpha \right) =\sum\displaylimits_{n=1}^{N}f_{n}\left(
x,y,A+\sigma \right) \Psi _{n}\left( \alpha \right) ,\left( x,y,\alpha
\right) \in \left( 0,1\right) ^{3},  \label{8}
\end{equation}%
\begin{equation}
F\left( x,y,A+\sigma \right) =\left( f_{1},f_{2},...,f_{N}\right) ^{T}\left(
x,y,A+\sigma \right) ,\left( x,y\right) \in \left( 0,1\right) ^{2},
\label{9}
\end{equation}%
Keeping in mind (\ref{4.18}) and (\ref{4.19}), we define the spaces $%
C_{N}^{1}\left( \overline{\Omega }\right) $, $C_{N,0}^{1}\left( \overline{%
\Omega }\right) $ of $N-$D vector functions $W\left( \mathbf{x}\right) $
defined in (\ref{4.19}) as 
\begin{equation*}
C_{x,y,N}^{1}\left( \overline{\Omega }\right) =\left\{ 
\begin{array}{c}
W\left( \mathbf{x}\right) :\left\Vert W\right\Vert _{C_{x,y,N}^{1}\left( 
\overline{\Omega }\right) } \\ 
=\max_{n\in \left[ 1,N\right] }\left( \left\Vert w_{n}\right\Vert _{C\left( 
\overline{\Omega }\right) }+\left\Vert w_{nx}\right\Vert _{C\left( \overline{%
\Omega }\right) }+\left\Vert w_{ny}\right\Vert _{C\left( \overline{\Omega }%
\right) }\right) <\infty%
\end{array}%
\right\} ,
\end{equation*}%
\begin{equation*}
C_{x,y,N,0}^{1}\left( \overline{\Omega }\right) =\left\{ W\in
C_{x,y,N}^{1}\left( \overline{\Omega }\right) :W\mid _{\Gamma }=W\left(
x,y,A\right) =0\right\} .
\end{equation*}

Keeping in mind (\ref{4.8})-(\ref{9}), substitute functions $w$, $g$ and
their corresponding first derivatives with respect to $x,y,\alpha $ in
equation (\ref{4.6}). Next, multiply the resulting equation sequentially by
functions $\Psi _{n}\left( \alpha \right) ,n=1,...,N$ and integrate with
respect to $\alpha \in \left( 0,1\right) $. Then multiply both sides of
obtained system of nonlinear integral differential equations by the matrix $%
M_{N}^{-1},$ where the matrix $M_{N}$ was introduced in section 3. We obtain%
\begin{equation}
W\left( \mathbf{x}\right) =M_{N}^{-1}P\left(
W,W_{x},W_{y},G,G_{x},G_{y},F_{x},F_{y},\mathbf{x}\right) ,\mathbf{x}\in
\Omega ,  \label{4.22}
\end{equation}%
\begin{equation}
W\mid _{\Gamma }=W\left( x,y,A\right) =0,  \label{4.23}
\end{equation}%
where $P$ is the $N-$D vector function, 
\begin{equation}
P=\left( P_{1},...,P_{N}\right) ^{T}\left(
W,W_{x},W_{y},G_{x},G_{y},F_{x},F_{y},\mathbf{x}\right) ,\mathbf{x}\in
\Omega ,  \label{4.24}
\end{equation}%
\begin{equation*}
P_{n}\left( W,W_{x},W_{y},G_{x},G_{y},F_{x},F_{y},\mathbf{x}\right) =
\end{equation*}%
\begin{equation}
=\int\displaylimits_{0}^{1}\Psi _{n}\left( \alpha \right) \frac{\partial }{%
\partial \alpha }\left( -\int\displaylimits_{z}^{A+\sigma }\frac{u_{0x}+w_{x}+g_{x}%
}{2\sqrt{u_{0}+w+g}}dt+f_{x}\left( x,y,A,\alpha \right) \right) ^{2}d\alpha
\label{4.240}
\end{equation}%
\begin{equation*}
+\int\displaylimits_{0}^{1}\Psi _{n}\left( \alpha \right) \frac{\partial }{%
\partial \alpha }\left( -\int\displaylimits_{z}^{A+\sigma }\frac{u_{0y}+w_{y}+g_{y}%
}{2\sqrt{u_{0}+w+g}}dt+f_{y}\left( x,y,A,\alpha \right) \right) ^{2}d\alpha .
\end{equation*}%
Thus, we have obtained the desired boundary value problem (\ref{4.22})-(\ref%
{4.240}) for the system of nonlinear coupled integral differential
equations. Below we focus on this problem.

\subsection{The positivity of the function $\left( u_{0}+w+g\right) \left( 
\mathbf{x}, \alpha \right) $}

\label{sec:5.3}

It follows from (\ref{4.6}) and (\ref{4.13}) that we need the function $%
\left( u_{0}+w+g\right) \left( \mathbf{x},\alpha \right) $ to be positive.
We discuss this issue in the current section.

Using (\ref{3.01}), (\ref{4.1}), (\ref{4.8}), (\ref{4.13}) and (\ref{4.19}),
define the set of functions $K$ as: 
\begin{equation}
K=\left\{ w\left( \mathbf{x},\alpha \right) :=\left( w+g\right) \left( 
\mathbf{x},\alpha \right) >0,\text{(\ref{4.16}) holds, \ }W\in
C_{x,y,N,0}^{1}\left( \overline{\Omega }\right) \right\} ,  \label{4.25}
\end{equation}%
where $\left( \mathbf{x},\alpha \right) \in \overline{\Omega }\times \left[
0,1\right] .$ Then by (\ref{4.8}) and (\ref{4.13}) 
\begin{equation}
\left( u_{0}+w+g\right) \left( \mathbf{x},\alpha \right) \geq \frac{A^{2}}{%
A^{2}+2},\forall w\in K,\left( \mathbf{x},\alpha \right) \in \overline{%
\Omega }\times \left[ 0,1\right] .  \label{4.26}
\end{equation}%
To obtain a sufficient condition guaranteeing (\ref{4.25}) in terms of the
vector function $W$, we formulate Lemma 5.1.

\textbf{Lemma 5.1}.\emph{\ Let (\ref{4.16}) and (\ref{4.17}) hold. Consider
the vector function }$v^{N}\left( \mathbf{x}\right) =\left(
w_{1}+g_{1},w_{2}+g_{2},...,w_{N}+g_{N}\right) ^{T}\left( \mathbf{x}\right)
. $\emph{\ Let }$X_{N}$\emph{\ be the }$N\times N$\emph{\ matrix of Lemma
3.1. Consider the vector function }$\widetilde{v}^{N}\left( \mathbf{x}%
\right) =X_{N}^{T}\cdot v^{N}\left( \mathbf{x}\right) .$\emph{\ Let }$%
\widetilde{v}^{N}\left( \mathbf{x}\right) =\left( v_{1},...,v_{N}\right)
^{T}\left( \mathbf{x}\right) .$ \emph{Suppose that all functions }$%
v_{n}\left( \mathbf{x}\right) >0$\emph{\ in }$\overline{\Omega }.$\emph{\
Then the function }$w\in K$\emph{\ and, therefore, (\ref{4.26}) holds}$.$%
\emph{\ Also, the set }$K$\emph{\ is convex}$.$

\textbf{Proof}. The first part of this lemma follows immediately from Lemma
3.1. We now prove the convexity of the set $K.$ Suppose that two functions $%
w^{\left( 1\right) },w^{\left( 2\right) }\in K.$ Let the number $\theta \in
\left( 0,1\right) .$ Then by (\ref{4.25}) 
\begin{equation*}
\theta w^{\left( 1\right) }\left( \mathbf{x},\alpha \right) >-\theta g\left( 
\mathbf{x},\alpha \right) ,\mathbf{x}\in \overline{\Omega },\alpha \in \left[
0,1\right] ,
\end{equation*}%
\begin{equation*}
\left( 1-\theta \right) w^{\left( 2\right) }\left( \mathbf{x},\alpha \right)
>-\left( 1-\theta \right) g\left( \mathbf{x},\alpha \right) ,\mathbf{x}\in 
\overline{\Omega },\alpha \in \left[ 0,1\right] .
\end{equation*}%
Summing up these two inequalities, we obtain%
\begin{equation*}
\theta w^{\left( 1\right) }\left( \mathbf{x},\alpha \right) +\left( 1-\theta
\right) w^{\left( 2\right) }\left( \mathbf{x},\alpha \right) +g\left( 
\mathbf{x},\alpha \right) >0,\mathbf{x}\in \overline{\Omega },\alpha \in %
\left[ 0,1\right] .\text{ }\square
\end{equation*}

\subsection{Applying the multidimensional analog of Taylor formula}

\label{sec:5.4}

We specify in this section how the classical multidimensional analog of
Taylor formula \cite{V}\emph{\ }can be applied to the right hand side of
equation (\ref{4.22}). Let $R>0$ be an arbitrary number. Denote%
\begin{equation}
K\left( R\right) =\left\{ W:w\in K,\left\Vert W\right\Vert
_{C_{x,y,N}^{1}\left( \overline{\Omega }\right) }<R\right\} .  \label{4.260}
\end{equation}%
It follows from Lemma 5.1 and (\ref{4.260}) that $K\left( R\right) $ is a
convex set.

\textbf{Lemma 5.2}. \emph{\ Let }$W^{\left( 1\right) }$\emph{,}$W^{\left(
2\right) }\in K\left( R\right) $\emph{,let }$G^{\left( 1\right) },G^{\left(
2\right) }\emph{be}$ $\emph{the}$ $\emph{vector}$ $\emph{functions}$ $\emph{(%
\ref{4.20})}$ \emph{and} $F^{\left( 1\right) },F^{\left( 2\right) }$ \emph{%
be the vector functions in (\ref{9}). Based on (\ref{4.16}) and (\ref{4.19}%
), denote }$\widetilde{w}\left( \mathbf{x},\alpha \right) =w^{\left(
1\right) }\left( \mathbf{x},\alpha \right) -w^{\left( 2\right) }\left( 
\mathbf{x},\alpha \right) $\emph{, }$\widetilde{w}_{n}\left( \mathbf{x}%
\right) =w_{n}^{\left( 1\right) }\left( \mathbf{x}\right) -w_{n}^{\left(
2\right) }\left( \mathbf{x}\right) .$ \emph{Similarly, denote }$\widetilde{g}%
\left( \mathbf{x},\alpha \right) =g^{\left( 1\right) }\left( \mathbf{x}%
,\alpha \right) -g^{\left( 2\right) }\left( \mathbf{x},\alpha \right) ,$ 
\emph{where} $g^{\left( k\right) }\left( \mathbf{x},\alpha \right) $\emph{\
corresponds to the vector function }$G^{\left( k\right) }\left( \mathbf{x}%
\right) ,k=1,2$ \emph{via (\ref{4.17}), (\ref{4.20}).} \emph{Also, denote }%
\begin{equation*}
\widetilde{W}=W^{\left( 2\right) }-W^{\left( 1\right) }=\left( \widetilde{w}%
_{1},...,\widetilde{w}_{N}\right) ,\text{ }\widetilde{G}=G^{\left( 1\right)
}-G^{\left( 2\right) },\text{ }\widetilde{F}=F^{\left( 2\right) }-F^{\left(
1\right) }.
\end{equation*}%
\emph{\ And let functions }$f^{\left( 1\right) }=f^{\left( 1\right) }\left(
x,y,A+\sigma \right) $ \emph{and} $\widetilde{f}=\widetilde{f}\left(
x,y,A+\sigma \right) $ \emph{correspond to the vector functions} $F^{\left(
1\right) }\left( \mathbf{x}\right) $ \emph{and} $\widetilde{F}\left( \mathbf{%
x}\right) $ \emph{respectively via (\ref{8}), (\ref{9}).} \emph{Then the
following form of the multidimensional Taylor formula is valid: } 
\begin{equation*}
M_{N}^{-1}P\left( W^{\left( 2\right) },W_{x}^{\left( 2\right)
},W_{y}^{\left( 2\right) },G_{x}^{\left( 2\right) },G_{y}^{\left( 2\right)
},F_{x}^{\left( 2\right) },F_{y}^{\left( 2\right) },\mathbf{x}\right)
\end{equation*}%
\begin{equation*}
-M_{N}^{-1}P\left( W^{\left( 1\right) },W_{x}^{\left( 1\right)
},W_{y}^{\left( 1\right) },G_{x}^{\left( 1\right) },G_{y}^{\left( 1\right)
},F_{x}^{\left( 1\right) },F_{y}^{\left( 1\right) }\mathbf{x}\right) =
\end{equation*}%
\begin{equation*}
=\int\displaylimits_{z}^{A+\sigma }T_{0}\left( W^{\left( 1\right) },W_{x}^{\left(
1\right) },W_{y}^{\left( 1\right) },G_{x}^{\left( 1\right) },G_{y}^{\left(
1\right) },F_{x}^{\left( 1\right) },F_{y}^{\left( 1\right) },x,y,t\right) 
\widetilde{W}\left( x,y,t\right) dt
\end{equation*}%
\begin{equation*}
+\int\displaylimits_{z}^{A+\sigma }T_{1}\left( W^{\left( 1\right) },W_{x}^{\left(
1\right) },W_{y}^{\left( 1\right) },G_{x}^{\left( 1\right) },G_{y}^{\left(
1\right) },F_{x}^{\left( 1\right) },F_{y}^{\left( 1\right) },x,y,t\right) 
\widetilde{W}_{x}\left( x,y,t\right) dt
\end{equation*}%
\begin{equation*}
+\int\displaylimits_{z}^{A+\sigma }T_{2}\left( W^{\left( 1\right) },W_{x}^{\left(
1\right) },W_{y}^{\left( 1\right) },G_{x}^{\left( 1\right) },G_{y}^{\left(
1\right) },F_{x}^{\left( 1\right) },F_{y}^{\left( 1\right) },x,y,t\right) 
\widetilde{W}_{y}\left( x,y,t\right) dt
\end{equation*}%
\begin{equation}
+\int\displaylimits_{z}^{A+\sigma }T_{3}\left( W^{\left( i\right) },W_{x}^{\left(
i\right) },W_{y}^{\left( i\right) },G_{x}^{\left( i\right) },G_{y}^{\left(
i\right) },F_{x}^{\left( i\right) },F_{y}^{\left( i\right) },\widetilde{W},%
\widetilde{W}_{x},\widetilde{W}_{y},x,y,t\right) dt  \label{4.261}
\end{equation}%
\begin{equation*}
+\int\displaylimits_{z}^{A+\sigma }S_{1}\left( W^{\left( 1\right) },W_{x}^{\left(
1\right) },W_{y}^{\left( 1\right) },G_{x}^{\left( 1\right) },G_{y}^{\left(
1\right) },F_{x}^{\left( 1\right) },F_{y}^{\left( 1\right) },x,y,t\right) 
\widetilde{G}_{x}\left( x,y,t\right) dt
\end{equation*}%
\begin{equation*}
+\int\displaylimits_{z}^{A+\sigma }S_{2}\left( W^{\left( 1\right) },W_{x}^{\left(
1\right) },W_{y}^{\left( 1\right) },G_{x}^{\left( 1\right) },G_{y}^{\left(
1\right) },F_{x}^{\left( 1\right) },F_{y}^{\left( 1\right) },x,y,t\right) 
\widetilde{G}_{y}\left( x,y,t\right) dt
\end{equation*}%
\begin{equation*}
+S_{3}\left( W^{\left( 1\right) },W_{x}^{\left( 1\right) },W_{y}^{\left(
1\right) },G_{x}^{\left( 1\right) },G_{y}^{\left( 1\right) },F_{x}^{\left(
1\right) },F_{y}^{\left( 1\right) },\mathbf{x}\right) \widetilde{F}%
_{x}\left( x,y\right)
\end{equation*}%
\begin{equation*}
+S_{4}\left( W^{\left( 1\right) },W_{x}^{\left( 1\right) },W_{y}^{\left(
1\right) },G_{x}^{\left( 1\right) },G_{y}^{\left( 1\right) },F_{x}^{\left(
1\right) },F_{y}^{\left( 1\right) }\mathbf{x}\right) \widetilde{F}_{y}\left(
x,y\right) ,
\end{equation*}%
\emph{where }$\mathbf{x}\in \Omega ,i=1,2$ \emph{in }$T_{3}$\emph{\ means
that this matrix depends on both vector functions }$\left( W^{\left(
1\right) },W_{x}^{\left( 1\right) },W_{y}^{\left( 1\right) },G_{x}^{\left(
1\right) },G_{y}^{\left( 1\right) },F_{x}^{\left( 1\right) },F_{y}^{\left(
1\right) }\right) $\emph{\ and }$\left( W^{\left( 2\right) },W_{x}^{\left(
2\right) },W_{y}^{\left( 2\right) },G_{x}^{\left( 2\right) },G_{y}^{\left(
2\right) },F_{x}^{\left( 2\right) },F_{y}^{\left( 2\right) }\right) .$ \emph{%
Also, all elements of matrices }$T_{j},S_{k},j=0,1,2,3;k=1,2,3,4$\emph{\ are
continuous functions of their variables.\ Furthermore, the following
estimates hold} 
\begin{equation}
\left\vert T_{j}\left( W^{\left( i\right) },W_{x}^{\left( i\right)
},W_{y}^{\left( i\right) },G_{x}^{\left( i\right) },G_{y}^{\left( i\right)
},F_{x}^{\left( i\right) },F_{y}^{\left( i\right) },x,y,t\right) \right\vert
\leq C_{1};j=0,1,2;\left( x,y,t\right) \in \overline{\Omega },  \label{4.262}
\end{equation}%
\begin{equation}
\left\vert T_{3}\left( W^{\left( i\right) },W_{x}^{\left( i\right)
},W_{y}^{\left( i\right) },G_{x}^{\left( i\right) },G_{y}^{\left( i\right)
},F_{x}^{\left( i\right) },F_{y}^{\left( i\right) },\widetilde{W},\widetilde{%
W}_{x},\widetilde{W}_{y},x,y,t\right) \right\vert  \label{4.263}
\end{equation}%
\begin{equation*}
\leq C_{1}\left( \widetilde{W}^{2}+\widetilde{W}_{x}^{2}+\widetilde{W}%
_{y}^{2}\right) ,\left( x,y,t\right) \in \overline{\Omega },
\end{equation*}%
\begin{equation}
\left\vert S_{k}\left( W^{\left( 1\right) },W_{x}^{\left( 1\right)
},W_{y}^{\left( 1\right) },G_{x}^{\left( 1\right) },G_{y}^{\left( 1\right)
},F_{x}^{\left( 1\right) },F_{y}^{\left( 1\right) },\mathbf{x}\right)
\right\vert \leq C_{1};k=1,...,4;\mathbf{x}\in \overline{\Omega },
\label{4.264}
\end{equation}%
\emph{where the number }%
\begin{equation}
C_{1}=C_{1}\left( N,R,\max_{i=1,2}\left\Vert G^{\left( i\right) }\right\Vert
_{C_{x,y,N}^{1}\left( \overline{\Omega }\right) },\max_{i=1,2}\left\Vert
F^{\left( i\right) }\right\Vert _{C_{x,y,N}^{1}\left( \overline{\Omega }%
\right) }\right) >0  \label{4.265}
\end{equation}%
\emph{\ depends only on listed parameters. Estimates }(\emph{\ref{4.262})-(%
\ref{4.265}) mean estimates for each element of corresponding matrices. In
terms of the integration with respect to }$\alpha ,$\emph{\ matrices }$%
T_{j},S_{j},j=0,1,2,3$ \emph{depend on the integrals of the form}%
\begin{equation*}
\int\displaylimits_{0}^{1}\Psi _{n}\left( \alpha \right) \mu _{jnkm}\left( \mathbf{%
x},\left( u_{0},u_{0x},u_{0y}\right) \left( x,y,t,\alpha \right) ,\Psi
_{k}\left( \alpha \right) ,\Psi _{m}\left( \alpha \right) ,\Psi _{k}^{\prime
}\left( \alpha \right) ,\Psi _{m}^{\prime }\left( \alpha \right) \right)
d\alpha ,\text{ }
\end{equation*}%
\emph{where }$n,k,m=1,...,N$ \emph{and} $\mu _{jnkm}$\emph{\ are continuous
functions of their variables. }

\textbf{Proof.} Below $C_{1}>0$ denotes different constants depending only
on parameters listed in (\ref{4.265}). For $\mathbf{x}\in \Omega $ and $%
\alpha \in \left( 0,1\right) $ consider the functions $\xi \left( \mathbf{x}%
,\alpha \right) $ and $\eta \left( \mathbf{x},\alpha \right) $ defined as 
\begin{equation}
\xi \left( w,\mathbf{x},\alpha \right) =\frac{u_{0x}+w_{x}+g_{x}}{2\sqrt{%
u_{0}+w+g}}\left( \mathbf{x},\alpha \right) ,\text{ }\eta \left( w,\mathbf{x}%
,\alpha \right) =\frac{u_{0y}+w_{y}+g_{y}}{2\sqrt{u_{0}+w+g}}\left( \mathbf{x%
},\alpha \right) ,  \label{4.27}
\end{equation}%
where functions $w,g$ have the forms (\ref{4.16}), (\ref{4.17}). Then $\xi
\left( w^{\left( 2\right) },\mathbf{x},\alpha \right) =\xi \left( w^{\left(
1\right) }+\widetilde{w},\mathbf{x},\alpha \right) $ and $\eta \left(
w^{\left( 2\right) },\mathbf{x},\alpha \right) =\eta \left( w^{\left(
1\right) }+\widetilde{w},\mathbf{x},\alpha \right) .$

The convexity of the set $K\left( R\right) $ allows us to use the
multidimensional analog of Taylor formula in the following form:%
\begin{equation*}
\xi \left( w^{\left( 2\right) },\mathbf{x},\alpha \right) =\xi \left(
w^{\left( 1\right) },\mathbf{x},\alpha \right) +s_{1}\left( w^{\left(
1\right) },\mathbf{x},\alpha \right) \widetilde{w}+s_{2}\left( w^{\left(
1\right) },\mathbf{x},\alpha \right) \widetilde{w}_{x}+
\end{equation*}%
\begin{equation}
s_{3}\left( w^{\left( 1\right) },\mathbf{x},\alpha \right) \widetilde{g}%
_{x}\left( \mathbf{x},\alpha \right) +s_{4}\left( w^{\left( 1\right)
},w^{\left( 2\right) },\mathbf{x},\alpha \right) \widetilde{w}^{2}
\label{4.28}
\end{equation}%
\begin{equation*}
+s_{5}\left( w^{\left( 1\right) },w^{\left( 2\right) },\mathbf{x},\alpha
\right) \widetilde{w}_{x}^{2}+s_{6}\left( w^{\left( 1\right) },w^{\left(
2\right) },\mathbf{x},\alpha \right) \widetilde{w}_{x}\widetilde{w},
\end{equation*}%
where $s_{j}$ are continuous functions of their listed variables.
Furthermore, 
\begin{equation}
\left\vert s_{i}\left( w^{\left( 1\right) },\mathbf{x},\alpha \right)
\right\vert ,\left\vert s_{j}\left( w^{\left( 1\right) },w^{\left( 2\right)
},\mathbf{x},\alpha \right) \right\vert \leq C_{1},\text{ where }i=1,2,3%
\text{ and }j=4,5,6.  \label{4.280}
\end{equation}%
Next, substituting (\ref{4.16}), (\ref{4.17}) and (\ref{4.20}) in (\ref{4.28}%
), we obtain%
\begin{equation*}
\xi \left( w^{\left( 2\right) },\mathbf{x},\alpha \right) =\xi
_{1}+\sum\displaylimits_{k=1}^{N}\left[ s_{1}\widetilde{w}_{k}\left( \mathbf{x}%
\right) +s_{2}\widetilde{w}_{kx}\left( \mathbf{x}\right) +s_{3}\widetilde{g}%
_{kx}\left( \mathbf{x}\right) \right] \Psi _{k}\left( \alpha \right)
\end{equation*}%
\begin{equation}
+\sum\displaylimits_{k,m=1}^{N}\left[ s_{4}\left( \widetilde{w}_{k}\widetilde{w}%
_{m}\right) \left( \mathbf{x}\right) +s_{5}\left( \widetilde{w}_{kx}%
\widetilde{w}_{mx}\right) \left( \mathbf{x}\right) +s_{6}\left( \widetilde{w}%
_{kx}\widetilde{w}_{m}\right) \left( \mathbf{x}\right) \right] \Psi
_{k}\left( \alpha \right) \Psi _{m}\left( \alpha \right) ,  \label{4.29}
\end{equation}%
where for brevity $\xi _{1}=\xi \left( w^{\left( 1\right) },\mathbf{x}%
,\alpha \right) .$ Hence, it follows from (\ref{4.27})-(\ref{4.29}) that the
second line of (\ref{4.240}) can be rewritten as 
\begin{equation*}
\int\displaylimits_{0}^{1}\Psi _{n}\left( \alpha \right) \frac{\partial }{\partial
\alpha }\left( -\int\displaylimits_{z}^{A+\sigma }\frac{u_{0x}+w_{x}+g_{x}}{2\sqrt{%
u_{0}+w+g}}dt+f_{x}\left( x,y,A,\alpha \right) \right) ^{2}d\alpha
=\int\displaylimits_{0}^{1}\Psi _{n}\left( \alpha \right) \times
\end{equation*}%
\begin{equation}
\times \frac{\partial }{\partial \alpha }\left( -\int\displaylimits_{z}^{A+\sigma
}\left( \xi _{1}+\sum\displaylimits_{k=1}^{N}\left(
V_{k}+\sum\displaylimits_{m=1}^{N}V_{km}\Psi _{m}\left( \alpha \right) \right)
\Psi _{k}\left( \alpha \right) \right) dt+\left( f_{x}^{\left( 1\right) }+%
\widetilde{f}_{x}\right) \right) ^{2}d\alpha ,  \label{4.30}
\end{equation}%
\begin{equation}
V_{k}\left( w^{\left( 1\right) },\mathbf{x},\alpha ,\widetilde{w}_{n}\right)
=s_{1}\widetilde{w}_{k}\left( \mathbf{x}\right) +s_{2}\widetilde{w}%
_{kx}\left( \mathbf{x}\right) +s_{3}\widetilde{g}_{kx}\left( \mathbf{x}%
\right) ,  \label{4.31}
\end{equation}%
\begin{equation}
V_{km}\left( w^{\left( 1\right) },w^{\left( 2\right) },\mathbf{x},\alpha
\right) =s_{4}\left( \widetilde{w}_{k}\widetilde{w}_{m}\right) \left( 
\mathbf{x}\right) +s_{5}\left( \widetilde{w}_{kx}\widetilde{w}_{mx}\right)
\left( \mathbf{x}\right) +s_{6}\left( \widetilde{w}_{kx}\widetilde{w}%
_{m}\right) \left( \mathbf{x}\right) .  \label{4.32}
\end{equation}%
Similar formulas are obviously valid for the third line of (\ref{4.240}).
Thus, formulas (\ref{4.240}), (\ref{4.27})-(\ref{4.32}) imply (\ref{4.261})-(%
\ref{4.265}). $\square $

In Lemma 5.2 we have used second order terms in the Taylor formula. In
addition, we will also need the formula which uses only linear terms. The
proof of Lemma 5.3 is completely similar to the proof of Lemma 5.2.

\textbf{Lemma 5.3.} \emph{Assume that conditions of Lemma 5.2 hold. Then,
the following analog of the final increment formula is valid }%
\begin{equation*}
M_{N}^{-1}P\left( W^{\left( 2\right) },W_{x}^{\left( 2\right)
},W_{y}^{\left( 2\right) },G_{x}^{\left( 2\right) },G_{y}^{\left( 2\right)
},F_{x}^{\left( 2\right) },F_{y}^{\left( 2\right) },\mathbf{x}\right)
\end{equation*}%
\begin{equation*}
-M_{N}^{-1}P\left( W^{\left( 1\right) },W_{x}^{\left( 1\right)
},W_{y}^{\left( 1\right) },G_{x}^{\left( 1\right) },G_{y}^{\left( 1\right)
},F_{x}^{\left( 1\right) },F_{y}^{\left( 1\right) }\mathbf{x}\right)
\end{equation*}%
\begin{equation*}
=\int\displaylimits_{z}^{A+\sigma }Y_{0}\left( W^{\left( i\right) },W_{x}^{\left(
i\right) },W_{y}^{\left( i\right) },G_{x}^{\left( i\right) },G_{y}^{\left(
i\right) },F_{x}^{\left( i\right) },F_{y}^{\left( i\right) },x,y,t\right) 
\widetilde{W}\left( x,y,t\right) dt
\end{equation*}%
\begin{equation*}
+\int\displaylimits_{z}^{A+\sigma }Y_{1}\left( W^{\left( i\right) },W_{x}^{\left(
i\right) },W_{y}^{\left( i\right) },G_{x}^{\left( i\right) },G_{y}^{\left(
i\right) },F_{x}^{\left( i\right) },F_{y}^{\left( i\right) },x,y,t\right) 
\widetilde{W}_{x}\left( x,y,t\right) dt
\end{equation*}%
\begin{equation}
+\int\displaylimits_{z}^{A+\sigma }Y_{2}\left( W^{\left( i\right) },W_{x}^{\left(
i\right) },W_{y}^{\left( i\right) },G_{x}^{\left( i\right) },G_{y}^{\left(
i\right) },F_{x}^{\left( i\right) },F_{y}^{\left( i\right) },x,y,t\right) 
\widetilde{W}_{y}\left( x,y,t\right) dt  \label{4.33}
\end{equation}%
\begin{equation*}
+\int\displaylimits_{z}^{A+\sigma }\widehat{S}_{1}\left( W^{\left( 1\right)
},W_{x}^{\left( 1\right) },W_{y}^{\left( 1\right) },G_{x}^{\left( 1\right)
},G_{y}^{\left( 1\right) },F_{x}^{\left( 1\right) },F_{y}^{\left( 1\right)
},x,y,t\right) \widetilde{G}_{x}\left( x,y,t\right) dt
\end{equation*}%
\begin{equation*}
+\int\displaylimits_{z}^{A+\sigma }\widehat{S}_{2}\left( W^{\left( 1\right)
},W_{x}^{\left( 1\right) },W_{y}^{\left( 1\right) },G_{x}^{\left( 1\right)
},G_{y}^{\left( 1\right) },F_{x}^{\left( 1\right) },F_{y}^{\left( 1\right)
},x,y,t\right) \widetilde{G}_{y}\left( x,y,t\right) dt
\end{equation*}%
\begin{equation*}
+\widehat{S}_{3}\left( W^{\left( 1\right) },W_{x}^{\left( 1\right)
},W_{y}^{\left( 1\right) },G_{x}^{\left( 1\right) },G_{y}^{\left( 1\right)
},F_{x}^{\left( 1\right) },F_{y}^{\left( 1\right) },\mathbf{x}\right) 
\widetilde{F}_{x}\left( x,y\right)
\end{equation*}%
\begin{equation*}
+\widehat{S}_{4}\left( W^{\left( 1\right) },W_{x}^{\left( 1\right)
},W_{y}^{\left( 1\right) },G_{x}^{\left( 1\right) },G_{y}^{\left( 1\right)
},F_{x}^{\left( 1\right) },F_{y}^{\left( 1\right) },\mathbf{x}\right) 
\widetilde{F}_{y}\left( x,y\right) ,\text{ }\mathbf{x}\in \Omega .
\end{equation*}%
\emph{where }$\mathbf{x}\in \Omega $, $i=1,2$\emph{,} \emph{all elements of
matrices }$Y_{j},j=0,1,2$\emph{\ are continuous functions of their variables
and the following estimates are valid for }$t\in \left[ z,A+\sigma \right]
,\left( x,y,t\right) ,\mathbf{x}\in \overline{\Omega }:$\emph{\ }%
\begin{equation*}
\left\vert \sum\displaylimits_{j=0}^{2}Y_{j}\left( W^{\left( i\right)
},W_{x}^{\left( i\right) },W_{y}^{\left( i\right) },G_{x}^{\left( i\right)
},G_{y}^{\left( i\right) },F_{x}^{\left( i\right) },F_{y}^{\left( i\right)
},x,y,t\right) \right\vert
\end{equation*}%
\begin{equation}
+\sum\displaylimits_{k=1}^{2}\left\vert \widehat{S}_{k}\left( W^{\left( 1\right)
},W_{x}^{\left( 1\right) },W_{y}^{\left( 1\right) },G_{x}^{\left( 1\right)
},G_{y}^{\left( 1\right) },F_{x}^{\left( 1\right) },F_{y}^{\left( 1\right)
},x,y,t\right) \right\vert  \label{4.34}
\end{equation}%
\begin{equation*}
+\sum\displaylimits_{k=3}^{4}\left\vert \widehat{S}_{k}\left( W^{\left( 1\right)
},W_{x}^{\left( 1\right) },W_{y}^{\left( 1\right) },G_{x}^{\left( 1\right)
},G_{y}^{\left( 1\right) },F_{x}^{\left( 1\right) },F_{y}^{\left( 1\right) },%
\mathbf{x}\right) \right\vert \leq C_{1}.
\end{equation*}

\section{Problem ( \ref{4.22}), ( \ref{4.23}) in the
Semidiscrete Form}

\label{sec:6}

We now rewrite equation (\ref{4.22}) in the form of finite differences with
respect to variables $x,y$ while keeping the continuous derivative with
respect to $z$. For brevity we keep the same grid step size $h>0$ in both
directions $x,y$. Consider partitions of the intervals $x\in \left(
0,1\right) ,y\in \left( 0,1\right) $ in small subintervals of the same
length $h$ with $B=1/h$ and corresponding semidiscrete sub-domains of the
domains $\Omega $ and $\overline{\Omega },$ 
\begin{equation*}
0=x_{0}<x_{1}<...<x_{B-1}<x_{B}=1,x_{i}-x_{i-1}=h,
\end{equation*}%
\begin{equation*}
0=y_{0}<y_{1}<...<y_{B-1}<y_{B}=1,y_{i}-y_{i-1}=h,
\end{equation*}%
\begin{equation}
\Omega ^{h}=\left\{ \mathbf{x}^{h}\left( z\right) =\left\{ \left(
x_{i},y_{j},z\right) \right\} _{i,j=1}^{B-1},z\in \left( A,A+\sigma \right)
\right\} .  \label{5.4}
\end{equation}%
\begin{equation}
\Omega _{1}^{h}=\left\{ \mathbf{x}^{h}\left( z\right) =\left\{ \left(
x_{i},y_{j},z\right) \right\} _{i,j=0}^{B},z\in \left( A,A+\sigma \right)
\right\} .  \label{5.5}
\end{equation}%
Hence, $\mathbf{x}^{h}\left( z\right) $ is a $z-$dependent vector. Its
dimension is $\left( B-2\right) ^{2}$ in the case of $\Omega ^{h}$ and $%
\left( B+1\right) ^{2}$ in the case of $\Omega _{1}^{h}.$ By (\ref{5.4})
only those points $\left( x_{i},y_{j},z\right) \in \Omega ^{h}$, which are
corresponding interior points of the domain $\Omega .$ On the other hand, in
addition to points of $\Omega ^{h},$ the semidiscrete domain $\Omega
_{1}^{h} $ contains boundary points which belong to the part $\Gamma $ of
the boundary $\partial \Omega .$ We assume below that 
\begin{equation}
h\in \left[ h_{0},1\right) ,h_{0}=const.\in \left( 0,1\right) .  \label{5.6}
\end{equation}

\textbf{Remark 6.1}. \emph{Unlike classical forward problems for PDEs, we do
not let the grid step size }$h$\emph{\ tend to zero. This is typical for
numerical methods for many inverse problems: due to their ill-posed nature,
see, e.g. \cite{KN,SKN}. In other words, the grid step size is often used as
the regularization parameter. }

Consider the $N-$D vector function $Q\left( \mathbf{x}\right) =\left(
Q_{1},....,Q_{N}\right) ^{T}\left( \mathbf{x}\right) $ with $Q_{n}\in
C\left( \overline{\Omega }\right) .$ We denote $Q^{h}\left( \mathbf{x}%
^{h}\left( z\right) \right) $ the trace of this vector function on the set $%
\overline{\Omega }_{1}^{h}$. Thus, 
\begin{equation}
Q^{h}\left( \mathbf{x}^{h}\left( z\right) \right) =\left(
Q_{1}^{h},...,Q_{N}^{h}\right) ^{T}\left( \mathbf{x}^{h}\left( z\right)
\right) =\left( \left( Q_{1}^{i,j}\left( z\right) \right)
_{i,j=0}^{B},...,\left( Q_{N}^{i,j}\left( z\right) \right)
_{i,j=0}^{B}\right) ^{T}  \label{5.60}
\end{equation}%
is the matrix depending on the variable $z\in \left[ A,A+\sigma \right] $.
Here $Q_{k}^{i,j}\left( z\right) =Q_{k}\left( x_{i},y_{j},z\right) ,$ where $%
k=1,...,N.$ Hence, by (\ref{4.16}), (\ref{4.17}), (\ref{8}) and (\ref{5.60})
the finite difference analogs of functions $w\left( \mathbf{x},\alpha
\right) $, $g\left( \mathbf{x},\alpha \right) ,f\left( x,y,A+\sigma \right) $
are: 
\begin{equation}
w^{h}\left( \mathbf{x}^{h}\left( z\right) ,\alpha \right)
=\sum\displaylimits_{n=1}^{N}w_{n}^{h}\left( \mathbf{x}^{h}\left( z\right) \right)
\Psi _{n}\left( \alpha \right) ,\text{ }\mathbf{x}^{h}\left( z\right) \in 
\overline{\Omega }_{1}^{h},\alpha \in \left[ 0,1\right] ,  \label{5.61}
\end{equation}%
\begin{equation}
g^{h}\left( \mathbf{x}^{h}\left( z\right) ,\alpha \right)
=\sum\displaylimits_{n=1}^{N}g_{n}^{h}\left( \mathbf{x}^{h}\left( z\right) \right)
\Psi _{n}\left( \alpha \right) ,\text{ }\mathbf{x}^{h}\left( z\right) \in 
\overline{\Omega }_{1}^{h},\alpha \in \left[ 0,1\right] .  \label{5.620}
\end{equation}%
\begin{equation}
f^{h}\left( \mathbf{x}^{h}\left( A+\sigma \right) ,\alpha \right)
=\sum\displaylimits_{n=1}^{N}f_{n}^{h}\left( \mathbf{x}^{h}\left( A+\sigma \right)
\right) \Psi _{n}\left( \alpha \right) ,\text{ }\mathbf{x}^{h}\left(
A+\sigma \right) \in \overline{\Omega }_{1}^{h},\alpha \in \left[ 0,1\right]
.  \label{5.63}
\end{equation}%
Next, by (\ref{4.19}), (\ref{4.20}), (\ref{9}) and (\ref{5.61})-(\ref{5.63})
the finite difference analogs of matrices $W$, $G$ and $F$ are%
\begin{equation}
W^{h}\left( \mathbf{x}^{h}\left( z\right) \right) =\left(
w_{1}^{h},...,w_{N}^{h}\right) ^{T}\left( \mathbf{x}^{h}\left( z\right)
\right) \text{, }\mathbf{x}^{h}\left( z\right) \in \overline{\Omega }%
_{1}^{h},  \label{5.64}
\end{equation}%
\begin{equation}
G^{h}\left( \mathbf{x}^{h}\left( z\right) \right) =\left(
g_{1}^{h},...,g_{N}^{h}\right) ^{T}\left( \mathbf{x}^{h}\left( z\right)
\right) ,\text{ }\mathbf{x}^{h}\left( z\right) \in \overline{\Omega }%
_{1}^{h},  \label{5.65}
\end{equation}%
\begin{equation}
F^{h}\left( \mathbf{x}^{h}\left( z\right) \right) =\left(
f_{1}^{h},f_{2}^{h},...,f_{N}^{h}\right) ^{T}\left( \mathbf{x}^{h}\left(
z\right) \right) ,z=A+\sigma ,\mathbf{x}^{h}\left( A+\sigma \right) \in 
\overline{\Omega }_{1}^{h}.  \label{5.66}
\end{equation}

For an arbitrary number $z\in \left[ A,A+\sigma \right] $ denote

$\Omega _{z}^{h}=\left\{ \mathbf{x}^{h}\left( z\right) =\left(
x_{i},y_{j},z\right) _{i,j=1}^{B-1},\text{ }z\text{ is fixed}\right\} .$ We
introduce semidiscrete functional spaces for matrices like $Q^{h}\left( 
\mathbf{x}\left( z\right) \right) ,$%
\begin{equation*}
C_{N}^{h}\left( \overline{\Omega }_{z}^{h}\right) =\left\{ Q^{h}\left( 
\mathbf{x}^{h}\left( z\right) \right) :\left\Vert Q^{h}\right\Vert
_{C_{N}^{h}\left( \overline{\Omega }_{z}^{h}\right) }\left( z\right)
=\max_{k\in \left[ 1,N\right] }\max_{i,j=1,...,B-1}\left\vert
Q_{k}^{i,j}\left( z\right) \right\vert <\infty \right\} ,
\end{equation*}%
\begin{equation*}
C_{N}^{h}\left( \overline{\Omega }^{h}\right) =\left\{ Q^{h}\left( \mathbf{x}%
^{h}\left( z\right) \right) :\left\Vert Q^{h}\right\Vert _{C_{N}^{h}\left( 
\overline{\Omega }^{h}\right) }=\max_{z\in \left[ A,A+\sigma \right]
}\left\Vert Q^{h}\right\Vert _{C_{N}^{h}\left( \overline{\Omega }%
_{z}^{h}\right) }\left( z\right) <\infty \right\} ,
\end{equation*}%
\begin{equation*}
C_{N}^{h}\left( \overline{\Omega }_{1,z}^{h}\right) =\left\{ Q^{h}\left( 
\mathbf{x}^{h}\left( z\right) \right) :\left\Vert Q^{h}\right\Vert
_{C^{h}\left( \overline{\Omega }_{1,z}^{h}\right) }\left( z\right)
=\max_{k\in \left[ 1,N\right] }\max_{i,j=0,...B}\left\vert Q_{k}^{i,j}\left(
z\right) \right\vert <\infty \right\} ,
\end{equation*}%
\begin{equation*}
C_{N}^{h}\left( \overline{\Omega }_{1}^{h}\right) =\left\{ Q^{h}\left( 
\mathbf{x}^{h}\left( z\right) \right) :\left\Vert Q^{h}\right\Vert
_{C_{N}^{h}\left( \overline{\Omega }_{1}^{h}\right) }=\max_{z\in \left[
A,A+\sigma \right] }\left\Vert Q^{h}\right\Vert _{C_{N}^{h}\left( \overline{%
\Omega }_{1z}^{h}\right) }\left( z\right) <\infty \right\} ,
\end{equation*}%
\begin{equation*}
C_{N,0}^{h}\left( \overline{\Omega }_{1}^{h}\right) =\left\{ Q^{h}\left( 
\mathbf{x}^{h}\left( z\right) \right) \in C_{N}^{h}\left( \overline{\Omega }%
_{1}^{h}\right) :Q^{h}\left( \mathbf{x}^{h}\left( z\right) \right) =0\text{
for }\mathbf{x}^{h}\left( z\right) \in \Gamma \cup \left\{ z=A\right\}
\right\} ,
\end{equation*}%
\begin{equation*}
L_{2,N}^{h}\left( \Omega _{z}^{h}\right) =\left\{ Q^{h}\left( \mathbf{x}%
^{h}\left( z\right) \right) :\left\Vert Q^{h}\right\Vert _{L_{2}^{h}\left(
\Omega _{z}^{h}\right) }^{2}\left( z\right)
=h^{2}\sum\displaylimits_{k=1}^{N}\sum\displaylimits_{i,j=1}^{B-1}\left(
Q_{k}^{i,j}\left( z\right) \right) ^{2}<\infty \right\} ,
\end{equation*}%
\begin{equation*}
L_{2,N}^{h}\left( \Omega ^{h}\right) =\left\{ Q^{h}\left( \mathbf{x}%
^{h}\left( z\right) \right) :\left\Vert Q^{h}\right\Vert _{L_{2,N}^{h}\left(
\Omega ^{h}\right) }^{2}=\int\displaylimits_{A}^{A+\sigma }\left\Vert
Q^{h}\right\Vert _{L_{2,N}^{h}\left( \Omega _{\theta }^{h}\right)
}^{2}\left( z\right) dz<\infty \right\} ,
\end{equation*}%
\begin{equation*}
H_{0,N}^{1,h}\left( \Omega _{1}^{h}\right) =\left\{ 
\begin{array}{c}
Q^{h}\left( \mathbf{x}^{h}\left( z\right) \right) :Q^{h}\left( \mathbf{x}%
^{h}\left( z\right) \right) \mid _{\Gamma \cup B_{A}}=0, \\ 
\left\Vert Q^{h}\right\Vert _{H_{0,N}^{1,h}\left( \Omega _{1}^{h}\right)
}^{2}= \\ 
=\int\displaylimits_{A}^{A+\sigma }\left[ \left\Vert Q^{h}\right\Vert
_{L_{2,N}^{h}\left( \Omega _{z}^{h}\right) }^{2}\left( z\right) +\left\Vert
\partial _{z}Q^{h}\right\Vert _{L_{2,N}^{h}\left( \Omega _{z}^{h}\right)
}^{2}\left( z\right) \right] dz<\infty%
\end{array}%
\right\} .
\end{equation*}

We approximate $x,y$ derivatives of the vector function $W\left( \mathbf{x}%
\right) $ via central finite differences \cite{Sam}. It is convenient to
write this in the equivalent form for the matrix function $W^{h}\left( 
\mathbf{x}\right) $ as%
\begin{equation}
W_{k,x}^{i,j}\left( z\right) =\frac{W_{k}^{i+1,j}\left( z\right)
-W_{k}^{i-1,j}\left( z\right) }{2h};i,j=1,...,B-1,  \label{5.8}
\end{equation}%
\begin{equation}
W_{k,y}^{i,j}\left( z\right) =\frac{W_{k}^{i,j+1}\left( z\right)
-W_{k}^{i,j-1}\left( z\right) }{2h},i,j=1,...,B-1,  \label{5.9}
\end{equation}%
\begin{equation}
W_{x}^{h}\left( \mathbf{x}^{h}\left( z\right) \right) =\left( \left(
W_{0,x}^{i,j}\left( z\right) \right) _{i,j=0}^{B},...,\left(
W_{0,x}^{i,j}\left( z\right) \right) _{i,j=0}^{B}\right) ^{T},  \label{5.10}
\end{equation}%
\begin{equation}
W_{y}^{h}\left( \mathbf{x}^{h}\left( z\right) \right) =\left( \left(
W_{0,y}^{i,j}\left( z\right) \right) _{i,j=0}^{B},...,\left(
W_{0,y}^{i,j}\left( z\right) \right) _{i,j=0}^{B}\right) ^{T}.  \label{5.11}
\end{equation}%
See (\ref{5.60}) for notations (\ref{5.10}) and (\ref{5.11}). Using (\ref%
{5.6}) and (\ref{5.8})-(\ref{5.11}), we obtain that there exists a constant $%
C_{2}=C_{2}\left( h_{0},N,\Omega ^{h}\right) >0$ depending only on listed
parameters such that 
\begin{equation}
\left\Vert Q_{x}^{h}\right\Vert _{C_{N}^{h}\left( \overline{\Omega }%
^{h}\right) },\left\Vert Q_{y}^{h}\right\Vert _{C_{N}^{h}\left( \overline{%
\Omega }^{h}\right) }\leq C_{2}\left\Vert Q^{h}\right\Vert _{C_{N}^{h}\left( 
\overline{\Omega }_{1}^{h}\right) },\text{ }\forall Q^{h}\in C_{N}^{h}\left( 
\overline{\Omega }_{1}^{h}\right) .  \label{5.12}
\end{equation}%
In addition, by embedding theorem $H_{0,N}^{1,h}\left( \Omega
_{1}^{h}\right) \subset C_{N,0}^{h}\left( \overline{\Omega }_{1}^{h}\right) $
and there exists a constant $C=C\left( A,\sigma ,h_{0},N\right) >0$ such that%
\begin{equation}
\left\Vert Q^{h}\right\Vert _{C_{N}^{h}\left( \overline{\Omega }%
_{1}^{h}\right) }\leq C\left\Vert Q^{h}\right\Vert _{H_{0,N}^{1,h}\left(
\Omega _{1}^{h}\right) },\text{ }\forall Q^{h}\in H_{0,N}^{1,h}\left( \Omega
_{1}^{h}\right) .  \label{5.121}
\end{equation}

Thus, using (\ref{5.60})-(\ref{5.66}), we obtain the following finite
difference analog of problem (\ref{4.22}), (\ref{4.23}):%
\begin{equation}
W^{h}\left( \mathbf{x}^{h}\left( z\right) \right) =M_{N}^{-1}P\left(
W^{h},W_{x}^{h},W_{y}^{h},G^{h},G_{x}^{h},G_{y}^{h},F^{h},F_{x}^{h},F_{y}^{h},%
\mathbf{x}^{h}\left( z\right) \right) ,\mathbf{x}^{h}\left( z\right) \in
\Omega ^{h},  \label{10}
\end{equation}%
\begin{equation}
W^{h}\left( \mathbf{x}^{h}\left( z\right) \right) \in H_{0,N}^{1,h}\left(
\Omega _{1}^{h}\right) .  \label{11}
\end{equation}%
Also, we assume that the vector functions $G^{h}\left( \mathbf{x}^{h}\left(
z\right) \right) $ and $F^{h}\left( \mathbf{x}^{h}\left( z\right) \right) $
in (\ref{5.65}) and (\ref{5.66}) are such that 
\begin{equation}
G^{h}\left( \mathbf{x}^{h}\left( z\right) \right) ,F^{h}\left( \mathbf{x}%
^{h}\left( z\right) \right) \in C_{N}^{h}\left( \overline{\Omega }%
_{1}^{h}\right) .  \label{13}
\end{equation}%
As above, let $R>0$ be an arbitrary number. We now introduce the finite
difference analogs of sets (\ref{4.25}) and (\ref{4.260}). Assume that (\ref%
{5.620}) holds. Then

\begin{equation}
K^{h}=\left\{ 
\begin{array}{c}
w^{h}=w^{h}\left( \mathbf{x}^{h}\left( z\right) ,\alpha \right) :=\left(
w^{h}+g^{h}\right) \left( \mathbf{x}^{h}\left( z\right) ,\alpha \right) >0,
\\ 
\forall \left( \mathbf{x}^{h}\left( z\right) ,\alpha \right) \in \overline{%
\Omega }_{1}^{h}\times \left[ 0,1\right] ,\text{ (\ref{5.61}) holds for }%
w^{h}\left( \mathbf{x}^{h}\left( z\right) ,\alpha \right) ,\text{ } \\ 
W^{h}\left( \mathbf{x}^{h}\left( z\right) \right) \in H_{0,N}^{1,h}\left(
\Omega _{1}^{h}\right)%
\end{array}%
\right\} ,  \label{14}
\end{equation}%
\begin{equation}
K^{h}\left( R\right) =\left\{ W^{h}\left( \mathbf{x}^{h}\left( z\right)
\right) :w^{h}\in K^{h},\left\Vert W^{h}\right\Vert _{H_{0,N}^{1,h}\left(
\Omega _{1}^{h}\right) }<R\right\} .  \label{15}
\end{equation}%
It follows from (\ref{5.121}), (\ref{14}) and (\ref{15}) that $K^{h}\left(
R\right) \subset K^{h}\subset C_{N,0}^{h}\left( \overline{\Omega }%
_{1}^{h}\right) .$

Similarly with Lemmata 3.1 and 5.1, Lemma 6.1 provides a sufficient
condition imposed on the components of the matrix $w^{h}\left( \mathbf{x}%
^{h}\left( z\right) ,\alpha \right) ,$ which guarantees that

$w^{h}\left( \mathbf{x}^{h}\left( z\right) ,\alpha \right) \in K^{h}.$

\textbf{Lemma 6.1}.\emph{\ Let the matrix }$w^{h}=w^{h}\left( \mathbf{x}%
^{h}\left( z\right) ,\alpha \right) \in K^{h}.$ Select a triple $\left(
i,j,z\right) $ with $i,j=0,...,B,z\in \left[ A,A+\sigma \right] $ and \emph{%
consider the vector }

$v\left( i,j,z\right) =\left( w_{1}+g_{1},w_{2}+g_{2},...,w_{N}+g_{N}\right)
^{T}\left( i,j,z\right) .$\emph{\ Let }$X_{N}$\emph{\ be the }$N\times N$%
\emph{\ matrix of Lemma 3.1. Consider the vector }$\widetilde{v}\left(
i,j,z\right) =X_{N}^{T}\cdot v\left( i,j,z\right) .$\emph{\ Let }

$\widetilde{v}\left( i,j,z\right) =\left( \widetilde{v}_{1},...,\widetilde{v}%
_{N}\right) ^{T}\left( i,j,z\right) .$ \emph{Suppose that all numbers }$%
\widetilde{v}_{n}\left( i,j,z\right) >0$ \emph{for all} $i,j=0,...,B,z\in %
\left[ A,A+\sigma \right] .$\emph{\ Then the function }$w^{h}\in K$\emph{\
and, therefore, by (\ref{4.8}) and (\ref{4.13}) the following analog of (\ref%
{4.26}) holds} for $\mathbf{x}^{h}\left( A+\sigma \right) ,\mathbf{x}%
^{h}\left( z\right) \in \overline{\Omega }_{1}^{h}:$ 
\begin{equation*}
u_{0}^{h}\left( \mathbf{x}\left( A+\sigma \right) \right) +w^{h}\left( 
\mathbf{x}^{h}\left( z\right) ,\alpha \right) +g^{h}\left( \mathbf{x}%
^{h}\left( z\right) ,\alpha \right) >\frac{A^{2}}{A^{2}+\_2}\text{ }.
\end{equation*}%
\emph{\ Also, the set }$K^{h}\left( R\right) $\emph{\ is convex}$.$

\textbf{Proof}. The first part of this lemma, the one about the positivity,
follows immediately from Lemma 3.1. Consider now two matrices $w^{1,h}\left( 
\mathbf{x}^{h}\left( z\right) ,\alpha \right) $,$w^{2,h}\left( \mathbf{x}%
^{h}\left( z\right) ,\alpha \right) \in K^{h}\left( R\right) .$ Let the
number $\theta \in \left[ 0,1\right] .$ Then one can prove completely
similarly with the proof of Lemma 5.1 that $\theta w^{1,h}\left( \mathbf{x}%
^{h}\left( z\right) ,\alpha \right) +\left( 1-\theta \right) w^{2,h}\left( 
\mathbf{x}^{h}\left( z\right) ,\alpha \right) \in K^{h}.$ Let $W^{1,h}\left( 
\mathbf{x}^{h}\left( z\right) \right) $ and $W^{2,h}\left( \mathbf{x}%
^{h}\left( z\right) \right) $ be two matrices corresponding to matrices $%
w^{1,h}\left( \mathbf{x}^{h}\left( z\right) ,\alpha \right) $ and

$w^{2,h}\left( \mathbf{x}^{h}\left( z\right) ,\alpha \right) $ respectively
via (\ref{5.64}).\emph{\ }The triangle inequality and (\ref{15}) imply that 
\begin{equation*}
\left\Vert \theta W^{1,h}\left( \mathbf{x}^{h}\left( z\right) \right)
+\left( 1-\theta \right) W^{2,h}\left( \mathbf{x}^{h}\left( z\right) ,\alpha
\right) \right\Vert _{H_{0,N}^{1,h}\left( \Omega _{1}^{h}\right) }
\end{equation*}%
\begin{equation*}
\leq \theta \left\Vert W^{1,h}\left( \mathbf{x}^{h}\left( z\right) \right)
\right\Vert _{H_{0,N}^{1,h}\left( \Omega _{1}^{h}\right) }+\left( 1-\theta
\right) \left\Vert W^{2,h}\left( \mathbf{x}^{h}\left( z\right) \right)
\right\Vert _{H_{0,N}^{1,h}\left( \Omega _{1}^{h}\right) }
\end{equation*}%
\begin{equation*}
<\theta R+\left( 1-\theta \right) R=R.\text{ }\square
\end{equation*}

Lemma 6.2 is a finite difference analog of Lemmata 5.2 and 5.3. The proof is
completely similar and is, therefore, omitted.

\textbf{Lemma 6.2}. \emph{Assume that (\ref{13}) holds. Then the direct
analogs of Lemmata 5.2 and 5.3, being applied to the right hand side of (\ref%
{10}), are true, provided that all functions involved in these lemmata are
replaced with their above semidiscrete analogs.} \emph{The constant }$C_{1}$%
\emph{\ in (\ref{4.265}) and (\ref{4.34}) should be replaced with the
constant }$\widetilde{C}_{1}$ \emph{depending only on listed parameters,
where} 
\begin{equation*}
\widetilde{C}_{1}=\widetilde{C}_{1}\left( h_{0},N,R,\max_{i=1,2}\left\Vert
G^{\left( i\right) }\right\Vert _{C_{x,y,N}\left( \overline{\Omega }\right)
},\max_{i=1,2}\left\Vert F^{\left( i\right) }\right\Vert
_{C_{x,y,N}^{1}\left( \overline{\Omega }\right) }\right) >0.
\end{equation*}

Suppose that we have found such a matrix $W^{h}\left( \mathbf{x}^{h}\left(
z\right) \right) \in K^{h}\left( R\right) ,$ that it solves problem (\ref{10}%
), (\ref{11}). Then, using (\ref{4.9}), (\ref{4.12}), (\ref{4.13}) and (\ref%
{5.61})-(\ref{5.63}), we set:%
\begin{equation}
w^{i,j}\left( z,\alpha \right) =\sum\displaylimits_{n=1}^{N}W_{n}^{i,j}\left(
z\right) \Psi _{n}\left( \alpha \right) ;i,j=0,...,B,z\in \left[ A,A+\sigma %
\right] ,\alpha \in \left( 0,1\right) ,  \label{5.15}
\end{equation}%
\begin{equation}
w^{h}\left( \mathbf{x}^{h}\left( z\right) ,\alpha \right) =\left(
w^{i,j}\left( z,\alpha \right) \right) _{i,j=0}^{B};\text{ }z\in \left[
A,A+\sigma \right] ,\alpha \in \left( 0,1\right) ,  \label{5.015}
\end{equation}%
\begin{equation}
v^{h}\left( \mathbf{x}^{h}\left( z\right) ,\alpha \right) =w^{h}\left( 
\mathbf{x}^{h}\left( z\right) ,\alpha \right) +g^{h}\left( \mathbf{x}%
^{h}\left( z\right) ,\alpha \right) ,\mathbf{x}^{h}\left( z\right) \in 
\overline{\Omega }_{1}^{h},\alpha \in \left( 0,1\right)  \label{5.115}
\end{equation}%
\begin{equation}
u^{h}\left( \mathbf{x}^{h}\left( z\right) ,\alpha \right) =u_{0}^{h}\left( 
\mathbf{x}^{h}\left( z\right) ,\alpha \right) +v^{h}\left( \mathbf{x}%
^{h}\left( z\right) ,\alpha \right) ,  \label{5.16}
\end{equation}%
$.$ Hence, by (\ref{4.8}), (\ref{14}) and (\ref{15}) 
\begin{equation}
u^{h}\left( \mathbf{x}^{h}\left( z\right) ,\alpha \right) >\frac{A^{2}}{%
A^{2}+2}>0,\forall \mathbf{x}^{h}\left( z\right) \in \overline{\Omega }%
_{1}^{h},\forall \alpha \in \left[ 0,1\right] .  \label{5.17}
\end{equation}%
Using (\ref{4.1}), (\ref{4.01}) and (\ref{5.17}), we set 
\begin{equation}
\tau _{z}^{h}\left( \mathbf{x}^{h}\left( z\right) ,\alpha \right) =\sqrt{%
u^{h}\left( \mathbf{x}^{h}\left( z\right) ,\alpha \right) }.  \label{5.18}
\end{equation}%
The semidiscrete analogs of formulas (\ref{4.3}) and (\ref{4.4}) are:%
\begin{equation}
\tau _{x}^{h}\left( \mathbf{x}^{h}\left( z\right) ,\alpha \right)
=-\int\displaylimits_{z}^{A+\sigma }\left( \frac{u_{x}^{h}}{2\sqrt{u^{h}}}\right)
\left( \mathbf{x}^{h}\left( t\right) ,\alpha \right) dt+f_{x}^{h}\left(
x,y,A+\sigma ,\alpha \right) ,\mathbf{x}^{h}\left( z\right) \in \Omega ,
\label{5.19}
\end{equation}%
\begin{equation}
\tau _{y}^{h}\left( \mathbf{x}^{h}\left( z\right) ,\alpha \right)
=-\int\displaylimits_{z}^{A+\sigma }\left( \frac{u_{y}^{h}}{2\sqrt{u^{h}}}\right)
\left( \mathbf{x}^{h}\left( t\right) ,\alpha \right) dt+f_{y}^{h}\left(
x,y,A+\sigma ,\alpha \right) ,\mathbf{x}^{h}\left( z\right) \in \Omega .
\label{5.20}
\end{equation}%
Next, using the original eikonal equation (\ref{2.9}), we set its
semidiscrete form as%
\begin{equation}
\left[ \left( \tau _{x}^{h}\right) ^{2}+\left( \tau _{y}^{h}\right)
^{2}+\left( \tau _{z}^{h}\right) ^{2}\right] \left( \mathbf{x}^{h}\left(
z\right) \right) =m^{h}\left( \mathbf{x}^{h}\left( z\right) \right) ,\text{ }%
\mathbf{x}^{h}\left( z\right) \in \overline{\Omega }^{h}.  \label{5.21}
\end{equation}

\textbf{Remark 6.2}. \emph{Equation (\ref{10}) with condition (\ref{11}) as
well conditions (\ref{5.6}), (\ref{5.61})-(\ref{5.63}) and the assumption
that the right hand side of (\ref{5.21}) is independent on the parameter }$%
\alpha $\emph{\ form our \textbf{approximate mathematical model} for the
TTTP formulated in section 2. One should expect that in practical
computations the left hand side of (\ref{5.21}) likely depends on }$\alpha .$
\emph{However, the numerical experience of \cite{KEIT,convnew,KN,SKN}, where
the basis }$\left\{ \Psi _{n}\left( \alpha \right) \right\} _{n=1}^{\infty }$%
\emph{\ was successfully used for four different inverse problems, shows
that, numerically, one should consider the average value of the left hand
side of (\ref{5.21}) with respect to }$\alpha .$

\textbf{Remarks 6.3:}

\textbf{(1)}\emph{\ It is well known that the problems like proving the
convergence of the numerical methods as ours when in (\ref{5.6}), (\ref{5.61}%
)-(\ref{5.63}) actual regularization parameters }$h_{0}\rightarrow 0$\emph{\
and }$N\rightarrow \infty $\emph{\ are, generally, extremely challenging
ones in the field of inverse problems. The fundamental underlying reason of
these challenges is the ill-posed nature of inverse problems. Therefore, we
do not analyze this type of convergence here.}

(\textbf{2}) \emph{Regardless on item 1, the author believes that the
success of numerical studies of \cite{KEIT,convnew,KN,SKN} indicates that
the numerical implementation of the method of this paper will likely lead to
good computational results. In this regard, we also mention here the work of
Guillement and Novikov} \emph{\cite{GN} for a linear inverse problem as as
well as the works of the group of Kabanikhin \cite{Kab1,Kab2,Kab3} for
nonlinear coefficient inverse problems. In all these four publications
truncated Fourier series were used to develop new numerical methods and to
test them numerically. Although convergence estimates at }$N\rightarrow
\infty $\emph{\ were not provided in \cite{GN,Kab1,Kab2,Kab3}, numerical
results are quite accurate ones.}

\section{Lipschitz Stability and Uniqueness}

\label{sec:7}

Based on (\ref{4.240}) as well as on (\ref{5.12}) and (\ref{5.121}), we
consider everywhere below the matrix equation (\ref{10}), as a system of
coupled nonlinear Volterra integral equations whose solution satisfies (\ref%
{11}). Denote%
\begin{equation}
\Phi ^{h}\left( \mathbf{x}^{h}\left( z\right) \right) =\left(
G^{h},G_{x}^{h},G_{y}^{h},F^{h},F_{x}^{h},F_{y}^{h}\right) \left( \mathbf{x}%
^{h}\left( z\right) \right) ;\text{ }\mathbf{x}^{h}\left( z\right) \in 
\overline{\Omega }^{h},\text{\ }G^{h},F^{h}\in C_{N}^{h}\left( \overline{%
\Omega }_{1}^{h}\right) .  \label{6.0}
\end{equation}

\textbf{Theorem 7.1} (Lipschitz stability and uniqueness). \emph{Let }$R>0$%
\emph{\ be an arbitrary number. Consider two matrices }$W^{\left( 1\right)
,h}\left( \mathbf{x}^{h}\left( z\right) \right) ,W^{\left( 2\right)
,h}\left( \mathbf{x}^{h}\left( z\right) \right) \in K^{h}\left( R\right) $%
\emph{. Suppose that these matrices generate two pairs of matrices in (\ref%
{5.65}), (\ref{5.66}) }$G^{\left( i\right) ,h}\left( \mathbf{x}^{h}\left(
z\right) \right) $ \emph{and} $F^{\left( i\right) ,h}\left( \mathbf{x}%
^{h}\left( A+\sigma \right) \right) ,i=1,2,$\emph{\ which satisfy conditions
(\ref{13}).Let }$\Phi ^{\left( i\right) ,h}\left( \mathbf{x}^{h}\left(
z\right) \right) ,i=1,2$ \emph{be the corresponding matrices as in (\ref{6.0}%
). Let} $m^{\left( 1\right) ,h}\left( \mathbf{x}\left( z\right) \right) $ 
\emph{and} $m^{\left( 2\right) ,h}\left( \mathbf{x}\left( z\right) \right) $ 
\emph{be the corresponding right hand sides of equality (\ref{5.21}). Assume
that functions }$m^{\left( 1\right) ,h}\left( \mathbf{x}\left( z\right)
\right) $ \emph{and} $m^{\left( 2\right) ,h}\left( \mathbf{x}\left( z\right)
\right) $ \emph{are independent on }$\alpha $ \emph{(Remarks 6.2). Then
there exists a constant}%
\begin{equation}
C_{3}=C_{3}\left( h_{0},A,\sigma ,R,N,\Omega ^{h},\left\Vert \Phi ^{\left(
1\right) ,h}\right\Vert _{C_{2N}^{h}\left( \overline{\Omega }_{1}^{h}\right)
},\left\Vert \Phi ^{\left( 2\right) ,h}\right\Vert _{C_{2N}^{h}\left( 
\overline{\Omega }_{1}^{h}\right) }\right) >0,  \label{6.1}
\end{equation}%
\emph{\ depending only on listed parameters such that the following
Lipschitz stability estimate holds} 
\begin{equation}
\left\Vert m^{\left( 1\right) ,h}-m^{\left( 2\right) ,h}\right\Vert
_{C_{1}^{h}\left( \overline{\Omega }^{h}\right) }\leq C_{3}\left\Vert \Phi
^{\left( 1\right) ,h}-\Phi ^{\left( 2\right) ,h}\right\Vert
_{C_{2N}^{h}\left( \overline{\Omega }_{1}^{h}\right) }.  \label{6.2}
\end{equation}%
\emph{In particular, if} $\Phi ^{\left( 1\right) ,h}\left( \mathbf{x}%
^{h}\left( z\right) \right) \equiv \Phi ^{\left( 2\right) ,h}\left( \mathbf{x%
}^{h}\left( z\right) \right) ,$ \emph{then} $m^{\left( 1\right) ,h}\left( 
\mathbf{x}^{h}\left( z\right) \right) \equiv m^{\left( 2\right) ,h}\left( 
\mathbf{x}^{h}\left( z\right) \right) ,$ \emph{i.e. uniqueness holds.}

\textbf{Proof}. In this proof, $C_{3}>0$ denotes different constants
depending on parameters listed in (\ref{6.1}). Denote 
\begin{equation}
\widetilde{W}^{h}=W^{\left( 2\right) ,h}-W^{\left( 1\right) ,h},\text{ }%
\widetilde{\Phi }^{h}=\Phi ^{\left( 2\right) ,h}-\Phi ^{\left( 1\right) ,h}.
\label{6.3}
\end{equation}%
It follows from Lemma 6.2, equality (\ref{4.33}), estimate (\ref{4.34}) of
Lemma 5.3, (\ref{5.12}), (\ref{10}), (\ref{11}), (\ref{5.15})-(\ref{5.17})
and (\ref{6.0}) that the following inequality with the Volterra integral
holds true 
\begin{equation*}
\left\vert \widetilde{W}^{h}\left( \mathbf{x}^{h}\left( z\right) \right)
\right\vert \leq C_{3}\int\displaylimits_{z}^{A+\sigma }\left\Vert \widetilde{W}%
^{h}\left( \mathbf{x}^{h}\left( t\right) \right) \right\Vert
_{C_{N}^{h}\left( \overline{\Omega }_{t}^{h}\right) }dt+C_{3}\left\Vert 
\widetilde{\Phi }^{h}\right\Vert _{C_{2N}^{h}\left( \overline{\Omega }%
_{1}^{h}\right) },
\end{equation*}%
where $\mathbf{x}\left( z\right) \in \overline{\Omega }_{1}^{h}$ and $z\in %
\left[ A,A+\sigma \right] .$ Hence, Gronwall's inequality leads to%
\begin{equation*}
\left\Vert \widetilde{W}^{h}\right\Vert _{C_{N}^{h}\left( \overline{\Omega }%
_{1}^{h}\right) }\leq C_{3}\left\Vert \widetilde{\Phi }^{h}\right\Vert
_{C_{2N}^{h}\left( \overline{\Omega }_{1}^{h}\right) }.
\end{equation*}%
This is the key estimate of this proof. Having this estimate, the target
estimate (\ref{6.2}) follows immediately from (\ref{5.64})-(\ref{5.66}), (%
\ref{5.12}), (\ref{5.15})-(\ref{5.21}), (\ref{6.0}) and (\ref{6.3}).
Uniqueness follows from (\ref{6.2}). $\square $

\section{Weighted Globally Strictly Convex Tikhonov-like Functional}

\label{sec:8}

\subsection{Estimating an integral}

\label{sec:8.1}

Let $\lambda >0$ be the parameter to be chosen later. We choose the
\textquotedblleft integral analog" of the CWF as 
\begin{equation}
\varphi _{\lambda }\left( z\right) =e^{2\lambda z}.  \label{9.1}
\end{equation}

\textbf{Lemma 8.1}. \emph{The following estimate holds for all }$\lambda >0$%
\emph{\ and for every function }$p\in L_{1}\left( A,A+\sigma \right) :$ 
\begin{equation*}
\int\displaylimits_{A}^{A+\sigma }e^{2\lambda z}\left( \int\displaylimits_{z}^{A+\sigma
}\left\vert p\left( y\right) \right\vert dy\right) dz\leq \frac{1}{2\lambda }%
\int\displaylimits_{A}^{A+\sigma }\left\vert p\left( z\right) \right\vert
e^{2\lambda z}dz.
\end{equation*}%
\textbf{Proof}. Interchanging the integrals, we obtain%
\begin{equation*}
I=\int\displaylimits_{A}^{A+\sigma }e^{2\lambda z}\left(
\int\displaylimits_{z}^{A+\sigma }\left\vert p\left( y\right) \right\vert
dy\right) dz=\int\displaylimits_{A}^{A+\sigma }\left\vert p\left( y\right)
\right\vert \left( \int\displaylimits_{A}^{y}e^{2\lambda z}dz\right) dy
\end{equation*}%
\begin{equation*}
=\frac{1}{2\lambda }\int\displaylimits_{A}^{A+\sigma }\left\vert p\left( y\right)
\right\vert \left( e^{2\lambda y}-e^{2\lambda A}\right) dy\leq \frac{1}{%
2\lambda }\int\displaylimits_{A}^{A+\sigma }\left\vert p\left( y\right)
\right\vert e^{2\lambda y}dy.\text{ \ }\square
\end{equation*}

\subsection{The functional}

\label{sec:8.2}

To solve problem (\ref{10}), (\ref{11}) numerically, we consider the
following minimization problem:

\textbf{Minimization Problem}. \emph{Fix an arbitrary number }$R>0$\emph{\
as well as the gird step size }$h\in \left[ h_{0},1\right) .$ \emph{Let }$%
\gamma >0$\emph{\ be the regularization parameter. Minimize the functional }$%
J_{\lambda ,\gamma }\left( W^{h}\right) $\emph{\ on the closed set }$%
\overline{K^{h}\left( R\right) },$\emph{\ where}%
\begin{equation}
J_{\lambda ,\gamma }\left( W^{h}\right) =e^{-2\lambda A}\left\Vert \left[
W^{h}\left( \mathbf{x}^{h}\left( z\right) \right) -M_{N}^{-1}P\left(
W^{h},W_{x}^{h},W_{y}^{h},\Phi ^{h},\mathbf{x}^{h}\left( z\right) \right) %
\right] e^{\lambda z}\right\Vert _{L_{2}^{h}\left( \Omega ^{h}\right) }^{2}
\label{9.4}
\end{equation}%
\begin{equation*}
+\gamma \left\Vert W^{h}\right\Vert _{H_{0}^{1,h}\left( \Omega
_{1}^{h}\right) }^{2}.
\end{equation*}

Here we took into account (\ref{9.1}). We use the multiplier $e^{-2\lambda
A} $ in order to balance two terms in the right hand side of (\ref{9.4}).
Note that since $R>0$ is an arbitrary number, then we do not impose a
smallness condition on the set $\overline{K^{h}\left( R\right) }$ where we
search for the solution of problem (\ref{10}), (\ref{11}). This is why we
are talking below about the \emph{global} strict convexity and the \emph{%
globally} convergent numerical method.

\textbf{Theorem 8.1} (global strict convexity).\emph{\ Let }$h_{0}$\emph{\
be the number defined in (\ref{5.6}) and let }$h\in \left[ h_{0},1\right) $%
\emph{. At every point }$W^{h}\in K^{h}\left( 2R\right) $ \emph{and for all }%
$\lambda >0,\gamma \in \left( 0,1\right) $ \emph{the functional }$J_{\lambda
,\gamma }\left( W^{h}\right) $\emph{\ has the Frechét derivative }$%
J_{\lambda ,\gamma }^{\prime }\left( W^{h}\right) \in H_{0}^{1,h}\left(
\Omega _{1}^{h}\right) $\emph{.\ Furthermore, this derivative is Lipschitz
continuous on }$K^{h}\left( 2R\right) $\emph{, i.e. there exists a constant }%
\begin{equation}
C_{4}=C_{4}\left( h_{0},A,\sigma ,R,N,\Omega ^{h},\left\Vert
G^{h}\right\Vert _{C_{N}^{h}\left( \overline{\Omega }_{1}^{h}\right)
},\left\Vert F^{h}\right\Vert _{C_{N}^{h}\left( \overline{\Omega }%
_{1}^{h}\right) }\right) >0  \label{9.5}
\end{equation}%
\emph{\ depending only on parameters listed in (\ref{9.5}) such that for all 
}$W^{\left( 1\right) ,h},W^{\left( 2\right) ,h}\in K^{h}\left( 2R\right) $%
\emph{\ }%
\begin{equation}
\left\Vert J_{\lambda ,\gamma }^{\prime }\left( W^{\left( 2\right)
,h}\right) -J_{\lambda ,\gamma }^{\prime }\left( W^{\left( 1\right)
,h}\right) \right\Vert _{H_{0}^{1,h}\left( \Omega _{1}^{h}\right) }\leq 
\overline{C}\left\Vert W^{\left( 2\right) ,h}-W^{\left( 1\right)
,h}\right\Vert _{H_{0}^{1,h}\left( \Omega _{1}^{h}\right) },  \label{9.6}
\end{equation}%
\emph{where the constant }$\overline{C}>0$\emph{\ depends on the same
parameters as ones listed in (\ref{9.5}) as well as on }$\lambda .$

\emph{In addition, there exists a sufficiently large number }$\lambda _{0}>1$%
\emph{\ depending on the same parameters as those listed in (\ref{9.5}) such
that for every }$\lambda \geq \lambda _{0}$\emph{\ and for every }$\gamma
\in \left( 0,1\right) $\emph{\ the functional }$J_{\lambda ,\gamma }\left(
W^{h}\right) $\emph{\ is strictly convex on the closed set }$\overline{%
K^{h}\left( R\right) }.$\emph{\ More precisely, the following estimate holds
for all }$W^{\left( 1\right) ,h},W^{\left( 2\right) ,h}\in \overline{%
K^{h}\left( R\right) }$\emph{\ }%
\begin{equation*}
J_{\lambda ,\gamma }\left( W^{\left( 2\right) ,h}\right) -J_{\lambda ,\gamma
}\left( W^{\left( 1\right) ,h}\right) -J_{\lambda ,\gamma }^{\prime }\left(
W^{\left( 1\right) ,h}\right) \left( W^{\left( 2\right) ,h}-W^{\left(
1\right) ,h}\right)
\end{equation*}%
\begin{equation}
\geq \frac{1}{8}\left\Vert W^{\left( 2\right) ,h}-W^{\left( 1\right)
,h}\right\Vert _{L_{2}^{h}\left( \Omega ^{h}\right) }^{2}+\gamma \left\Vert
W^{\left( 2\right) ,h}-W^{\left( 1\right) ,h}\right\Vert _{H_{0}^{1,h}\left(
\Omega _{1}^{h}\right) }^{2}.  \label{9.7}
\end{equation}

Below $C_{4}$ denotes different constants depending on parameters listed in (%
\ref{9.5})\emph{. }Theorem 8.2 follows immediately from (\ref{9.6}), (\ref%
{9.7}) and Lemma 2.1 of \cite{BakKlib}

\textbf{Theorem 8.2}. \emph{Suppose that conditions of Theorem 8.1 are in
place. Then for every }$\lambda \geq \lambda _{0}$\emph{\ and for every }$%
\gamma \in \left( 0,1\right) $\emph{\ there exists a single minimizer }$%
W_{\min }^{h}\in \overline{K^{h}\left( R\right) }$\emph{\ of the functional }%
$J_{\lambda ,\gamma }\left( W^{h}\right) $\emph{\ on the set }$\overline{%
K^{h}\left( R\right) }.$\emph{\ Furthermore,} 
\begin{equation}
J_{\lambda ,\gamma }^{\prime }\left( W_{\min }^{h}\right) \left(
W^{h}-W_{\min }^{h}\right) \geq 0,\text{ }\forall W^{h}\in \overline{%
K^{h}\left( R\right) }.  \label{9.8}
\end{equation}

Let $P_{\overline{K^{h}\left( R\right) }}:H_{0}^{1,h}\left( \Omega
_{1}^{h}\right) \rightarrow \overline{K^{h}\left( R\right) }$ be the
projection operator of the space $H_{0}^{1,h}\left( \Omega _{1}^{h}\right) $
on the closed set $\overline{K^{h}\left( R\right) }\subset H_{0}^{1,h}\left(
\Omega _{1}^{h}\right) .$ Consider an arbitrary point $W^{0,h}\in
K^{h}\left( R\right) .$ And minimize the functional $J_{\lambda ,\gamma
}\left( W^{h}\right) $ by the gradient projection method, which starts its
iterations at the point $W^{0,h},$ 
\begin{equation}
W_{n}^{h}=P_{\overline{K^{h}\left( R\right) }}\left( W_{n-1}^{h}-\rho
J_{\lambda ,\alpha }^{\prime }\left( W_{n-1}^{h}\right) \right) ,n=1,2,...
\label{9.9}
\end{equation}%
Here the number $\rho \in \left( 0,1\right) $. Theorem 8.3 follows from a
combination of Theorems 8.1 and 8.2 with Theorem 2.1 of \cite{BakKlib}.

\textbf{Theorem 8.3. }\emph{Suppose that conditions of Theorem 8.1 are in
place, }$\lambda \geq \lambda _{0}$\emph{\ and }$\gamma \in \left(
0,1\right) .$\emph{\ Let }$W_{\min }^{h}\in \overline{K^{h}\left( R\right) }$
\emph{be the minimizer listed in Theorem 8.2. Then there exists a
sufficiently small number }$\rho _{0}\in \left( 0,1\right) $\emph{\
depending on the same parameters as ones listed in (\ref{9.5}) such that for
every }$\rho \in \left( 0,\rho _{0}\right) $\emph{\ there exists a number }$%
\eta =\eta \left( \rho \right) \in \left( 0,1\right) $\emph{\ such that the
sequence (\ref{9.9}) converges to }$W_{\min }^{h},$%
\begin{equation}
\left\Vert W_{\min }^{h}-W_{n}^{h}\right\Vert _{H_{0,N}^{1,h}\left( \Omega
_{1}^{h}\right) }\leq \eta ^{n}\left\Vert W_{\min }^{h}-W_{0}^{h}\right\Vert
_{H_{0,N}^{1,h}\left( \Omega _{1}^{h}\right) }.  \label{9.10}
\end{equation}

In the regularization theory, the minimizer $W_{\min }^{h}$ of functional (%
\ref{9.4}) is called \textquotedblleft regularized solution", see, e.g. \cite%
{BK,T}. We now need to show that regularized solutions converge to the exact
one when the noise in the data tends to zero. Following the regularization
theory, we assume that there exists an exact, i.e. idealized, solution $%
W^{\ast ,h}\in K^{h}\left( R\right) $ of problem (\ref{10}), (\ref{11}) with
the noiseless data $\Phi ^{\ast ,h}\left( \mathbf{x}\left( z\right) \right)
, $ 
\begin{equation}
\Phi ^{\ast ,h}\left( \mathbf{x}\left( z\right) \right) =\left( G^{\ast
,h},G_{x}^{\ast ,h},G_{y}^{\ast ,h},F^{\ast ,h},F_{x}^{\ast ,h},F_{y}^{\ast
,h}\right) \left( \mathbf{x}^{h}\left( z\right) \right) ;\text{ }G^{\ast
,h},F^{\ast ,h}\in C_{N}^{h}\left( \overline{\Omega }_{1}^{h}\right) ,
\label{9.11}
\end{equation}%
where $\mathbf{x}^{h}\left( z\right) \in \overline{\Omega }^{h},$ see (\ref%
{6.0}). We assume that there exists the exact, i.e. idealized function $%
m^{\ast ,h}\left( \mathbf{x}^{h}\left( z\right) \right) ,$ $\mathbf{x}%
^{h}\left( z\right) \in \overline{\Omega }^{h}$ in (\ref{5.21}), which
produces the data (\ref{9.11}).

Let the number $\delta \in \left( 0,1\right) $ be the level of the error in
the data $G^{h},F^{h},$ i.e. 
\begin{equation}
\left\Vert G^{\ast ,h}-G^{h}\right\Vert _{C_{N}^{h}\left( \overline{\Omega }%
_{1}^{h}\right) },\left\Vert F^{\ast ,h}-F^{h}\right\Vert _{C_{N}^{h}\left( 
\overline{\Omega }_{1}^{h}\right) }<\delta .  \label{9.12}
\end{equation}%
Denote $\widetilde{G}^{h}=G^{\ast ,h}-G^{h},\widetilde{F}^{h}=F^{\ast
,h}-F^{h}.$ Then (\ref{5.12}) and (\ref{9.12}) imply that with a constant $%
C_{2}=C_{2}\left( h_{0},N,\Omega ^{h}\right) >0$ depending only on listed
parameters the following inequalities hold:%
\begin{equation}
\left\Vert \widetilde{G}^{h}\right\Vert _{C_{N}^{h}\left( \overline{\Omega }%
_{1}^{h}\right) },\left\Vert \widetilde{G}_{x}^{h}\right\Vert
_{C_{N}^{h}\left( \overline{\Omega }^{h}\right) },\left\Vert \widetilde{G}%
_{y}^{h}\right\Vert _{C_{N}^{h}\left( \overline{\Omega }^{h}\right)
}<C_{2}\delta ,  \label{9.13}
\end{equation}%
\begin{equation}
\left\Vert \widetilde{F}^{h}\right\Vert _{C_{N}^{h}\left( \overline{\Omega }%
_{1}^{h}\right) },\left\Vert \widetilde{F}_{x}^{h}\right\Vert
_{C_{N}^{h}\left( \overline{\Omega }^{h}\right) },\left\Vert \widetilde{F}%
_{y}^{h}\right\Vert _{C_{N}^{h}\left( \overline{\Omega }^{h}\right)
}<C_{2}\delta .  \label{9.14}
\end{equation}

Since $\delta \in \left( 0,1\right) $ and (\ref{9.12})-(\ref{9.14}) hold,
then, using the triangle inequality, we replace below the dependence of the
constant $C_{4}>0$ on $\left\Vert \Phi ^{h}\right\Vert _{C_{2N}^{h}\left( 
\overline{\Omega }_{1}^{h}\right) }$ in (\ref{9.5}) with the dependence of $%
C_{4}$ on $\left\Vert \Phi ^{\ast ,h}\right\Vert _{C_{2N}^{h}\left( 
\overline{\Omega }_{1}^{h}\right) }.$ Thus, everywhere below $C_{4}>0$
denotes different constants depending on the same parameters as ones listed
in (\ref{9.5}), except that $\left\Vert \Phi ^{h}\right\Vert
_{C_{2N}^{h}\left( \overline{\Omega }_{1}^{h}\right) }$ is replaced with $%
\left\Vert \Phi ^{\ast ,h}\right\Vert _{C_{2N}^{h}\left( \overline{\Omega }%
_{1}^{h}\right) }.$

Consider now the right hand side of equation (\ref{10}) for the case when $%
W^{h}$ is replaced with $W^{\ast ,h},$ whereas other arguments remain the
same. By Lemma 6.2, we can use finite difference analogs of (\ref{4.33}) and
(\ref{4.34}). In addition, we use now (\ref{9.11})-(\ref{9.14}). Thus, we
obtain 
\begin{equation*}
M_{N}^{-1}P\left( W^{\ast ,h},W_{x}^{\ast ,h},W_{y}^{\ast h},\Phi ^{h},%
\mathbf{x}^{h}\left( z\right) \right)
\end{equation*}%
\begin{equation}
=M_{N}^{-1}P\left( W^{\ast ,h},W_{x}^{\ast ,h},W_{y}^{\ast h},\Phi ^{\ast
,h},\mathbf{x}^{h}\left( z\right) \right) +\widehat{P}\left( \mathbf{x}%
^{h}\left( z\right) \right) ,\mathbf{x}^{h}\left( z\right) \in \overline{%
\Omega }^{h},  \label{9.15}
\end{equation}%
\begin{equation}
\left\Vert \widehat{P}\right\Vert _{C_{N}^{h}\left( \overline{\Omega }%
^{h}\right) }\leq C_{4}\delta .  \label{9.16}
\end{equation}

\textbf{Theorem 8.4}. \emph{Suppose that conditions of Theorem 8.1 are in
place and also that there exists an exact solution} $W^{\ast ,h}\in
K^{h}\left( R\right) $ o\emph{f problem (\ref{10}), (\ref{11}) with the
noiseless data }$\Phi ^{\ast ,h}$ \emph{as in (\ref{9.11}).} \emph{Let the }$%
\lambda _{0}>1$\emph{\ be the number chosen in Theorem 8.1. Fix an arbitrary
number }$\lambda =\lambda _{1}\geq \lambda _{0}$ \emph{in the functional }$%
J_{\lambda ,\gamma }\left( W^{h}\right) =J_{\lambda _{1},\gamma }\left(
W^{h}\right) .$ \emph{Just as in the regularization theory, set }$\gamma
=\gamma \left( \delta \right) =\delta ^{2}.$\emph{\ Let the numbers }$\rho
\in \left( 0,\rho _{0}\right) ,\rho _{0},\eta =\eta \left( \rho \right) \in
\left( 0,1\right) $\emph{\ be the same as in Theorem 8.3. Then the following
estimates are valid:}%
\begin{equation}
\left\Vert W^{h,\ast }-W_{\min }^{h}\right\Vert _{L_{2,N}^{h}\left( \Omega
^{h}\right) }\leq C_{4}\delta ,  \label{9.161}
\end{equation}%
\begin{equation}
\left\Vert W^{\ast ,h}-W_{n}^{h}\right\Vert _{L_{2,N}^{h}\left( \Omega
^{h}\right) }\leq C_{4}\delta +\eta ^{n}\left\Vert W_{\min
}^{h}-W_{0}^{h}\right\Vert _{H_{0,N}^{1,h}\left( \Omega _{1}^{h}\right) }.
\label{9.162}
\end{equation}%
\begin{equation}
\left\Vert m^{\ast ,h}-m_{n}^{h}\right\Vert _{L_{2,1}^{h}\left( \Omega
^{h}\right) }\leq C_{4}\delta +\eta ^{n}\left\Vert W_{\min
}^{h}-W_{0}^{h}\right\Vert _{H_{0,N}^{1,h}\left( \Omega _{1}^{h}\right) },
\label{9.163}
\end{equation}%
\emph{where the functions} $m_{n}^{h}\left( \mathbf{x}^{h}\left( z\right)
\right) $ \emph{are constructed from matrices }$W_{n}^{h}$\emph{\ via the
procedure described in section 6 with the final formula (\ref{5.21}).}

It follows from the above that we need to prove only Theorems 8.1 and 8.4.

\section{Proofs of Theorems 8.1 and 8.4}

\label{sec:9}

\subsection{Proof of Theorem 8.1}

\label{sec:9.1}

Denote $\left( ,\right) $ and $\left[ ,\right] $ the scalar products in the
spaces $L_{2,N}^{h}\left( \Omega ^{h}\right) $ and $H_{0,N}^{1,h}\left(
\Omega _{1}^{h}\right) $ respectively. Let $W^{\left( 1\right) ,h}$ and $%
W^{\left( 2\right) ,h}$ be two arbitrary points of the set $K^{h}\left(
R\right) .$ As above, denote $\widetilde{W}^{h}=W^{\left( 2\right)
,h}-W^{\left( 1\right) ,h}.$ Hence, $W^{\left( 2\right) ,h}=W^{\left(
1\right) ,h}+\widetilde{W}^{h}.$ Also, by (\ref{15}) 
\begin{equation}
\left\Vert \widetilde{W}^{h}\right\Vert _{H_{0,N}^{1,h}\left( \Omega
_{1}^{h}\right) }<R.  \label{8.1}
\end{equation}%
Hence, by (\ref{5.121}) and (\ref{8.1}) 
\begin{equation}
\left\Vert \widetilde{W}^{h}\right\Vert _{C_{N}^{h}\left( \overline{\Omega }%
^{h}\right) }\leq 2CR.  \label{8.2}
\end{equation}%
It follows from Lemma 6.2 and (\ref{4.261}) that 
\begin{equation*}
M_{N}^{-1}P\left( W^{\left( 2\right) ,h},W_{x}^{\left( 2\right)
,h},W_{y}^{\left( 2\right) ,h},\Phi ^{h},\mathbf{x}^{h}\left( z\right)
\right)
\end{equation*}%
\begin{equation*}
=M_{N}^{-1}P\left( W^{\left( 1\right) ,h},W_{x}^{\left( 1\right)
,h},W_{y}^{\left( 1\right) ,h},\Phi ^{h},\mathbf{x}^{h}\left( z\right)
\right)
\end{equation*}%
\begin{equation*}
+\int\displaylimits_{z}^{A+\sigma }T_{0}\left( W^{\left( 1\right)
,h},W_{x}^{\left( 1\right) ,h},W_{y}^{\left( 1\right) ,h},\Phi
^{h},x,y,t\right) \widetilde{W}^{h}\left( \mathbf{x}^{h}\left( t\right)
\right) dt
\end{equation*}%
\begin{equation}
+\int\displaylimits_{z}^{A+\sigma }T_{1}\left( W^{\left( 1\right)
,h},W_{x}^{\left( 1\right) ,h},W_{y}^{\left( 1\right) ,h},\Phi
^{h},x,y,t\right) \widetilde{W}_{x}^{h}\left( \mathbf{x}^{h}\left( t\right)
\right) dt  \label{8.3}
\end{equation}%
\begin{equation*}
+\int\displaylimits_{z}^{A+\sigma }T_{2}\left( W^{\left( 1\right)
,h},W_{x}^{\left( 1\right) ,h},W_{y}^{\left( 1\right) ,h},\Phi
^{h},x,y,t\right) \widetilde{W}_{y}^{h}\left( \mathbf{x}^{h}\left( t\right)
\right) dt
\end{equation*}%
\begin{equation*}
+\int\displaylimits_{z}^{A+\sigma }T_{3}\left( W^{\left( i\right)
,h},W_{x}^{\left( i\right) ,h},W_{y}^{\left( i\right) ,h},\Phi ^{h},%
\widetilde{W}^{h},\widetilde{W}_{x}^{h},\widetilde{W}_{y}^{h},\mathbf{x}%
^{h}\left( t\right) \right) dt,\text{ }\mathbf{x}^{h}\left( z\right) \in 
\overline{\Omega }^{h}.
\end{equation*}%
where $i=1,2$ and $T_{j},j=0,1,2,3$ are continuous functions of their
variables for $W^{\left( i\right) ,h}\in \overline{D^{h}\left( R\right) }.$
In addition, by (\ref{4.263}) and (\ref{8.2}) 
\begin{equation*}
\int\displaylimits_{z}^{A+\sigma }\left\vert T_{3}\left( W^{\left( i\right)
,h},W_{x}^{\left( i\right) ,h},W_{y}^{\left( i\right) ,h},\Phi ^{h},%
\widetilde{W}^{h},\widetilde{W}_{x}^{h},\widetilde{W}_{y}^{h},\mathbf{x}%
^{h}\left( t\right) \right) \right\vert
\end{equation*}%
\begin{equation}
\leq C_{4}\int\displaylimits_{z}^{A+\sigma }\left( \left( \widetilde{W}^{h}\right)
^{2}+\left( \widetilde{W}_{x}^{h}\right) ^{2}+\left( \widetilde{W}%
_{y}^{h}\right) ^{2}\right) \left( \mathbf{x}^{h}\left( t\right) \right) dt,%
\text{ }\mathbf{x}^{h}\left( z\right) \in \overline{\Omega }^{h}.
\label{8.4}
\end{equation}

Denote 
\begin{equation}
D_{1}^{h}\left( \mathbf{x}^{h}\left( z\right) \right) =W^{\left( 1\right)
,h}\left( \mathbf{x}^{h}\left( z\right) \right) -M_{N}^{-1}P\left( W^{\left(
1\right) ,h},W_{x}^{\left( 1\right) ,h},W_{y}^{\left( 1\right) ,h},\Phi ^{h},%
\mathbf{x}^{h}\left( z\right) \right) ,  \label{8.5}
\end{equation}%
\begin{equation}
D_{2}^{h}\left( \widetilde{W}^{h},\widetilde{W}_{x}^{h},\widetilde{W}%
_{y}^{h},\mathbf{x}^{h}\left( z\right) \right) =\text{the sum of lines
number 3,4,5 of (\ref{8.3}),}  \label{8.6}
\end{equation}%
\begin{equation}
D_{3}^{h}\left( \widetilde{W}^{h},\widetilde{W}_{x}^{h},\widetilde{W}%
_{y}^{h},\mathbf{x}^{h}\left( z\right) \right) =\text{the line number 6 of (%
\ref{8.3}).}  \label{8.7}
\end{equation}%
Then it follows from (\ref{9.4}) and (\ref{8.3})-(\ref{8.7}) that%
\begin{equation*}
J_{\lambda ,\gamma }\left( W^{\left( 2\right) ,h}\right) -J_{\lambda ,\gamma
}\left( W^{\left( 1\right) ,h}\right)
\end{equation*}%
\begin{equation*}
=2e^{-2\lambda \left( A+\sigma \right) }\left( D_{1}^{h},\widetilde{W}%
^{h}-D_{2}^{h}\left( \widetilde{W}^{h},\widetilde{W}_{x}^{h},\widetilde{W}%
_{y}^{h}\right) e^{2\lambda z}\right) +2\gamma \left[ \widetilde{W}%
^{h},W^{\left( 1\right) ,h}\right]
\end{equation*}%
\begin{equation}
-2e^{-2\lambda A}\left( D_{1}^{h},D_{3}^{h}\left( \widetilde{W}^{h},%
\widetilde{W}_{x}^{h},\widetilde{W}_{y}^{h}\right) e^{2\lambda z}\right)
\label{8.8}
\end{equation}%
\begin{equation*}
+e^{-2\lambda A}\left\Vert \left[ \widetilde{W}^{h}-D_{2}^{h}\left( 
\widetilde{W}^{h},\widetilde{W}_{x}^{h},\widetilde{W}_{y}^{h}\right)
-D_{3}^{h}\left( \widetilde{W}^{h},\widetilde{W}_{x}^{h},\widetilde{W}%
_{y}^{h}\right) \right] e^{\lambda z}\right\Vert _{L_{2,N}^{h}\left( \Omega
^{h}\right) }^{2}
\end{equation*}%
\begin{equation*}
+\gamma \left\Vert \widetilde{W}^{h}\right\Vert _{H_{0,N}^{1,h}\left( \Omega
_{1}^{h}\right) }^{2}.
\end{equation*}%
In addition, by Lemma 6.2, (\ref{4.261}), (\ref{4.262}), (\ref{5.121}),
Lemma 8.1, (\ref{8.2}) and (\ref{8.4})-(\ref{8.7}) the following estimates
hold:%
\begin{equation*}
\left\vert -2\left( D_{1}^{h},D_{3}^{h}\left( \widetilde{W}^{h},\widetilde{W}%
_{x}^{h},\widetilde{W}_{y}^{h}\right) e^{2\lambda z}\right) \right\vert
e^{-2\lambda A}
\end{equation*}%
\begin{equation}
\leq C_{4}e^{-2\lambda A}\int\displaylimits_{A}^{A+\sigma }\left(
\int\displaylimits_{z}^{A+\sigma }\left\Vert \widetilde{W}^{h}\right\Vert
_{L_{2.N}\left( \Omega _{t}^{h}\right) }^{2}dt\right) e^{2\lambda z}dz\leq 
\frac{C_{4}}{\lambda }e^{-2\lambda A}\left\Vert \widetilde{W}^{h}e^{\lambda
z}\right\Vert _{L_{2.N}\left( \Omega ^{h}\right) }^{2}.  \label{8.9}
\end{equation}%
Similarly using Cauchy-Schwarz inequality, Lemma 8.1 and assuming that the
parameter $\lambda \geq \lambda _{0},$ where $\lambda _{0}$ is a
sufficiently large number depending on the same parameters as ones listed in
(\ref{9.5}), we estimate from the below the sum of the fourth and fifth
lines of (\ref{8.8}) as:%
\begin{equation*}
e^{-2\lambda A}\left\Vert \left[ \widetilde{W}^{h}-D_{2}^{h}\left( 
\widetilde{W}^{h},\widetilde{W}_{x}^{h},\widetilde{W}_{y}^{h}\right)
-V_{3}^{h}\left( \widetilde{W}^{h},\widetilde{W}_{x}^{h},\widetilde{W}%
_{y}^{h}\right) \right] e^{\lambda z}\right\Vert _{L_{2,N}^{h}\left( \Omega
^{h}\right) }^{2}
\end{equation*}%
\begin{equation*}
+\gamma \left\Vert \widetilde{W}^{h}\right\Vert _{H_{0,N}^{1,h}\left( \Omega
_{1}^{h}\right) }^{2}
\end{equation*}%
\begin{equation}
\geq \frac{1}{2}e^{-2\lambda A}\left\Vert \widetilde{W}^{h}e^{\lambda
z}\right\Vert _{L_{2.N}\left( \Omega ^{h}\right) }^{2}-\frac{C_{4}}{\lambda }%
e^{-2\lambda A}\left\Vert \widetilde{W}^{h}e^{\lambda z}\right\Vert
_{L_{2.N}\left( \Omega ^{h}\right) }^{2}+\gamma \left\Vert \widetilde{W}%
^{h}\right\Vert _{H_{0,N}^{1,h}\left( \Omega _{1}^{h}\right) }^{2}
\label{8.10}
\end{equation}%
\begin{equation*}
\geq \frac{1}{4}e^{-2\lambda A}\left\Vert \widetilde{W}^{h}e^{\lambda
z}\right\Vert _{L_{2.N}\left( \Omega ^{h}\right) }^{2}+\gamma \left\Vert 
\widetilde{W}^{h}\right\Vert _{H_{0,N}^{1,h}\left( \Omega _{1}^{h}\right)
}^{2}.
\end{equation*}

It follows from (\ref{8.3}), (\ref{8.5}) and (\ref{8.6}) that the expression
in the second line of (\ref{8.8}) is linear with respect to $\widetilde{W}%
^{h},$%
\begin{equation}
Lin\left( \widetilde{W}^{h}\right) =2e^{-2\lambda A}\left( D_{1}^{h},%
\widetilde{W}^{h}-D_{2}^{h}\left( \widetilde{W}^{h},\widetilde{W}_{x}^{h},%
\widetilde{W}_{y}^{h}\right) e^{2\lambda z}\right) +2\gamma \left[ 
\widetilde{W}^{h},W^{\left( 1\right) ,h}\right] .  \label{8.11}
\end{equation}%
Obviously $\left\vert Lin\left( \widetilde{W}^{h}\right) \right\vert \leq
C_{4}\left\Vert \widetilde{W}^{h}\right\Vert _{H_{0,N}^{1,h}\left( \Omega
_{1}^{h}\right) }.$ Hence, $Lin\left( \widetilde{W}^{h}\right)
:H_{0,N}^{1,h}\left( \Omega _{1}^{h}\right) \rightarrow \mathbb{R}$ is a
bounded linear functional. Hence by Riesz theorem there exists a matrix $%
\Theta \in H_{0,N}^{1,h}\left( \Omega _{1}^{h}\right) $ such that%
\begin{equation}
Lin\left( \widetilde{W}^{h}\right) =\left[ \Theta ,\widetilde{W}^{h}\right] ,%
\text{ }\forall \widetilde{W}^{h}\in H_{0,N}^{1,h}\left( \Omega
_{1}^{h}\right) .  \label{8.12}
\end{equation}%
Besides, it follows from the above that%
\begin{equation*}
J_{\lambda ,\gamma }\left( W^{\left( 1\right) ,h}+\widetilde{W}^{h}\right)
-J_{\lambda ,\gamma }\left( W^{\left( 1\right) ,h}\right) -\left[ \Theta ,%
\widetilde{W}^{h}\right] =o\left( \left\Vert \widetilde{W}^{h}\right\Vert
_{H_{0,N}^{1,h}\left( \Omega _{1}^{h}\right) }\right) ,
\end{equation*}
as $\left\Vert \widetilde{W}^{h}\right\Vert _{H_{0,N}^{1,h}\left( \Omega
_{1}^{h}\right) }\rightarrow 0.$ Hence, $\Theta =J_{\lambda ,\gamma
}^{\prime }\left( W^{\left( 1\right) ,h}\right) \in H_{0,N}^{1,h}\left(
\Omega _{1}^{h}\right) $ is the Frechét derivative of the functional $%
J_{\lambda ,\gamma }$ at the point $W^{\left( 1\right) ,h}\in K^{h}\left(
R\right) .$ The existence of the Frechét derivative in the larger set $%
K^{h}\left( 2R\right) $ can be proven completely similarly. The Lipschitz
continuity property (\ref{9.6}) of the Frechét derivative $J_{\lambda
,\gamma }^{\prime }$ can be proven similarly with the proof of Theorem 3.1
of \cite{BakKlib}. Therefore,\ we omit this proof for brevity.

Furthermore, using (\ref{8.8})-(\ref{8.12}) and recalling that $\widetilde{W}%
^{h}=W^{\left( 2\right) ,h}-W^{\left( 1\right) ,h},$ we obtain for
sufficiently large $\lambda _{0}>1$ and for $\lambda \geq \lambda _{0}$ 
\begin{equation*}
J_{\lambda ,\gamma }\left( W^{\left( 2\right) ,h}\right) -J_{\lambda ,\gamma
}\left( W^{\left( 1\right) ,h}\right) -J_{\lambda ,\gamma }^{\prime }\left(
W^{\left( 1\right) ,h}\right) \left( W^{\left( 2\right) ,h}-W^{\left(
1\right) ,h}\right)
\end{equation*}%
\begin{equation*}
\geq \frac{1}{8}\left\Vert \widetilde{W}^{h}\right\Vert _{L_{2.N}\left(
\Omega ^{h}\right) }^{2}+\gamma \left\Vert \widetilde{W}^{h}\right\Vert
_{H_{0,N}^{1,h}\left( \Omega _{1}^{h}\right) }^{2},
\end{equation*}%
which is the same as the target estimate (\ref{9.7}) of this theorem. $%
\square $

\subsection{Proof of Theorem 8.4}

\label{sec:9.2}

Recall that in this theorem we fix and arbitrary number $\lambda =\lambda
_{1}\geq \lambda _{0},$ where $\lambda _{0}$ is the number of Theorem 8.1.
To indicate the dependence of the functional $J_{\lambda _{1},\gamma }$ on $%
\Phi ^{h},$ we write in this proof $J_{\lambda _{1},\gamma }\left(
W^{h},\Phi ^{h}\right) .$

First, we consider $J_{\lambda _{1},\gamma }\left( W^{\ast ,h},\Phi ^{\ast
,h}\right) .$ Since 
\begin{equation}
W^{\ast ,h}\left( \mathbf{x}^{h}\left( z\right) \right) -M_{N}^{-1}P\left(
W^{\ast ,h},W_{x}^{\ast ,h},W_{y}^{\ast h},\Phi ^{\ast ,h},\mathbf{x}%
^{h}\left( z\right) \right) =0,\text{ }\mathbf{x}^{h}\left( z\right) \in
\Omega ^{h},  \label{8.13}
\end{equation}%
then, using (\ref{9.4}), we obtain%
\begin{equation}
J_{\lambda _{0},\gamma }\left( W^{\ast ,h},\Phi ^{\ast ,h}\right) =\gamma
\left\Vert W^{\ast ,h}\right\Vert _{H_{0}^{1,h}\left( \Omega _{1}^{h}\right)
}^{2}\leq \gamma R^{2}.  \label{8.14}
\end{equation}%
Next, by (\ref{9.4}), (\ref{9.15}), (\ref{9.16}), (\ref{8.13}) and (\ref%
{8.14})%
\begin{equation}
J_{\lambda _{1},\gamma }\left( W^{\ast ,h},\Phi ^{h}\right) =J_{\lambda
_{1},\gamma }\left( W^{\ast ,h},\Phi ^{\ast h}\right) +Z\left( W^{\ast
,h},W^{h},\Phi ^{\ast h},\Phi ^{h}\right) ,  \label{8.15}
\end{equation}%
where $Z\left( W^{\ast ,h},W^{h},\Phi ^{\ast h},\Phi ^{h}\right) $ satisfies
the following estimate 
\begin{equation}
\left\vert Z\left( W^{\ast ,h},W^{h},\Phi ^{\ast h},\Phi ^{h}\right)
\right\vert \leq C_{4}\delta ^{2}e^{2\lambda _{1}\sigma }.  \label{8.16}
\end{equation}%
Since $\lambda _{1}$ is a fixed arbitrary number such that $\lambda _{1}\geq
\lambda _{0}$ and since $\lambda _{0}$ depends on the same parameters as
those listed in (\ref{9.5}) for $C_{4},$ then recalling that $\gamma =\delta
^{2},$ we obtain from (\ref{8.14})-(\ref{8.16})%
\begin{equation}
J_{\lambda _{1},\gamma }\left( W^{\ast ,h},\Phi ^{h}\right) \leq C_{4}\delta
^{2}.  \label{8.17}
\end{equation}

Next, (\ref{9.7}) implies that 
\begin{equation}
J_{\lambda _{1},\gamma }\left( W^{\ast ,h},\Phi ^{h}\right) -J_{\lambda
_{1},\gamma }\left( W_{\min }^{h},\Phi ^{h}\right) -J_{\lambda _{1},\gamma
}^{\prime }\left( W_{\min }^{h},\Phi ^{h}\right) \left( W^{\ast ,h}-W_{\min
}^{h}\right)  \label{8.18}
\end{equation}%
\begin{equation*}
\geq \frac{1}{8}\left\Vert W^{\ast ,h}-W_{\min }^{h}\right\Vert
_{L_{2}^{h}\left( \Omega ^{h}\right) }^{2}.
\end{equation*}%
Since by (\ref{9.8}) $-J_{\lambda _{0},\gamma }^{\prime }\left( W_{\min
}^{h},\Phi ^{h}\right) \left( W^{\ast ,h}-W_{\min }^{h}\right) \leq 0,$
then, using (\ref{8.17}) and (\ref{8.18}), we obtain%
\begin{equation}
\left\Vert W^{\ast ,h}-W_{\min }^{h}\right\Vert _{L_{2}^{h}\left( \Omega
^{h}\right) }^{2}\leq C_{4}\delta ^{2}.  \label{8.19}
\end{equation}%
The first target estimate (\ref{9.161}) of this theorem follows from (\ref%
{8.19}). Combining (\ref{9.10}) and (\ref{9.161}) with the procedure of
section 6, which led to (\ref{5.21}), we obtain two other target estimates (%
\ref{9.162}) and (\ref{9.163}) of Theorem 8.4. $\square $

\textbf{Acknowledgments.} This work was supported by US Army Research
Laboratory and US Army Research Office grant W911NF-19-1-0044. The author is
grateful to Professor Vladimir G. Romanov for many very fruitful discussions.

%\bibliographystyle{siamplain}
%\bibliography{references}

\end{document}